\documentclass[11pt]{article}
\title{Global divergence of spatial coalescents}

\author{Omer Angel\thanks{University of British Columbia; supported in part
    by NSERC.} %
  \and{Nathana\"el Berestycki\thanks{Cambridge University; supported in
      part by EPSRC grant EP/G055068/1.}} %
  \and{Vlada Limic\thanks{CNRS, UMR 6632; supported in part by Alfred
      P.~Sloan Research Fellowship, and in part by ANR MAEV grant.}}}

\date{November 2009}

\usepackage{amsmath,amsthm,amssymb,amsfonts,mathdots}
\usepackage[margin=2.5cm]{geometry}
\usepackage[ansinew]{inputenc}
\usepackage[dvips]{graphics,color}
\usepackage{graphicx}
\usepackage{hyperref}

\newtheorem{thm}{Theorem}[section]
\newtheorem{lemma}[thm]{Lemma}

\theoremstyle{definition}

\newtheorem*{rem*}{Remark}

\newcommand{\thmref}[1]{Theorem~\ref{T:#1}}
\newcommand{\lemref}[1]{Lemma~\ref{L:#1}}

\newcommand{\eps}{\varepsilon}
\newcommand{\E}{\mathbb{E}}
\renewcommand{\P}{\mathbb{P}}
\newcommand{\R}{\mathbb{R}}
\newcommand{\Z}{\mathbb{Z}}
\newcommand{\N}{\mathbb{N}}
\newcommand{\slt}{\preceq}
\newcommand{\sgt}{\succeq}
\newcommand{\FF}{\mathcal{F}}
\newcommand{\GG}{\mathcal{G}}
\newcommand{\cP}{\mathcal{P}}

\newcommand{\La}{\Lambda}

\newcommand{\tN}{N^*}

\def\crw{^{\diamond}}

\DeclareMathOperator{\Var}{Var}
\DeclareMathOperator{\Poi}{Poisson}
\DeclareMathOperator{\Vol}{Vol}

\DeclareMathOperator{\Bin}{Bin}
\DeclareMathOperator{\Tow}{Tow}
\DeclareMathOperator{\Geom}{Geom}

\newcommand{\indic}[1]{\mathbf{1}_{\{#1\}}}
\newcommand{\indica}[1]{\mathbf{1}_{#1}}

\begin{document}

\maketitle

\begin{abstract}
We study several fundamental properties of a class of stochastic processes
called spatial $\Lambda$-coalescents.
In these models, a number of particles perform
independent random walks on some underlying graph $G$. In addition,
particles on the same vertex merge randomly according to a given
coalescing mechanism. A remarkable property of mean-field coalescent
processes is that they may come
down from infinity, meaning that, starting with an infinite number of
particles, only a finite number remains after any positive amount of
time, almost surely.
We show here however that, in the spatial setting, on any infinite and
bounded-degree graph, the total number of particles will always remain
infinite at all times,
almost surely. Moreover, if $G=\Z^d$, and the coalescing
mechanism is Kingman's coalescent, then starting with $N$ particles at
the origin, the total number of particles remaining is of
order $(\log^* N)^d$ at any fixed positive time (where $\log^*$ is the
inverse tower function). At sufficiently large times the total number
of particles is of order $(\log^* N)^{d-2}$, when $d>2$.
We provide parallel results in the recurrent case $d=2$.
The spatial Beta-coalescents behave similarly,
where $\log\log N$ is replacing $\log^* N$.
\end{abstract}

\setcounter{tocdepth}{1}
\tableofcontents
\newpage

\section{Introduction}

\subsection{Motivation and main results}

The theory of stochastic coalescent processes has expanded considerably in
the last decade, as a consequence of their deep connections to
population genetics, spin glass models and polymers.
In theoretical population genetics, coalescents arise as natural models of
merging of ancestral lineages (see, for example, \cite{durrettbookDNA,
  bertoin}). A particular $\La$-coalescent, usually called the
Bolthausen-Sznitman coalescent, is thought to be an important object for
describing the conjectured universal ultrametric structure of numerous
mean-field spin glass models including the Sherrington-Kirkpatrick model
(see \cite{boltsznit, bovkur, talagrand}). The same coalescent has also
been recently linked in \cite{derridaetal} to scaling limits of directed
polymers.

The $\La$-coalescents are stochastic processes taking values in
$\cP$, the space of partitions of $\N=\{1,2,\ldots\}$.
In the current context, each class (or block) in
$\Pi \in \cP$ can and will be thought of as a particle.
For any coalescent process $(\Pi_{t},t\geq 0)$ it is true that
$\Pi_{t+s}$ is a coarsening of $\Pi_t$, for any $t,s>0$. There is a natural
semi-group structure $(\cP,\star)$, where $\Pi\star\Pi'$ is the result of
merging the blocks of $\Pi$ ``according to'' the partition
$\Pi'$. The $\La$-coalescents can be canonically characterized as the L\'evy
processes in $(\cP,\star)$, with the property that
no two coagulation events occur simultaneously. The above L\'evy property
corresponds to the fact that
$\Pi_{t+s} = \Pi_t \star \Pi'_s$,
for all $t,s\ge 0$, where $(\Pi'_t, t\ge0)$ is an independent, identically
distributed process (see, e.g.\ \cite[\S 3.1.3]{coalnotes} for
details).
A direct construction of $\La$-coalescents from \cite{pit99} is now considered
standard in the probability literature. As will be discussed in more detail
below, each $\La$-coalescent corresponds uniquely to a finite measures
$\La$ on $[0,1]$: for instance, the case where $\La = \delta_{\{0\}}$ gives
the well-known \emph{Kingman} coalescent from mathematical population
genetics, while the uniform measure on $[0,1]$ gives the
Bolthausen-Sznitman coalescent. These two processes, as well as the
more general Beta-coalescents are especially interesting
and amenable for analysis, since they are self-similar in a
certain sense made precise by the results of \cite{bbs2}.

The present work is devoted to the study of several fundamental properties
of a more general class of models, introduced in \cite{LS} and called the
spatial $\La$-coalescents. In this setting, particles (i.e., partition
classes) are positioned on some underlying locally finite graph $G=(V,E)$.
The dynamics of the process is enriched in the presence of the geographical
structure in two ways. First, the particles move as independent continuous
time simple random walks on $G$. Second, stochastically independent
$\La$-coalescence takes place on each site of $G$. More precisely, at any
given time, only particles that are on a same site can coagulate. Moreover,
at every site the coalescence mechanism is that of the original
(mean-field) $\La$-coalescent.
The spatial Kingman coalescent is a natural model of an interacting
particle system where particles perform independent random walks on an
underlying graph, and any pair of particles coalesce at rate $1$, as long
as they are located at the same site.

The spatial $\La$-coalescent processes are particularly well suited to
model merging of ancestral lineages for a population that is evolving in a
geographical space $G$, where the spatial motion of individuals is taken
into account. In this way, geographical factors such as isolation and
overpopulation can influence the dynamics, making it a more realistic model
for long-term population behavior. In the above interpretation, the
vertices of the graph are referred to as {\em demes} and represent a
discretization of physical space. Each edge represents potential migratory
routes between two adjacent demes.

While the mean-field $\La$-coalescent processes are relatively well
understood at this point, even basic properties of their spatial
counterparts are much more delicate to analyze. Intuitively, the difficulty
comes from the fact that the two ingredients in the dynamics, the
coalescence and the migration, affect the particles in the system in
opposite directions: the spatial motion makes particles diffuse away from
one another, and the coalescence keeps them together. Indeed, our results
show that the competition between these two forces can be very tight. Our
main results provide information about the limiting behavior of spatial
$\La$-coalescents as the initial number of particles $n$ tends to infinity,
at both small and large time-scales. We consider the case where initially
all the particles are located at the origin $o$ of $G$. For some of our
results, the only assumptions on $G$ are that it is connected and has
bounded degree $\Delta = \max_{v \in V}{\rm degree}(v)<\infty$. However,
several of our more precise results on the asymptotic behavior are restricted
to the setting where $G$ is the $d$-dimensional lattice $\Z^d$.

Define the function $\log^* n$ as the inverse $\log^* n := \inf\{m\ge 1:
\text{Tow}(m,1)\ge n\}$ of the tower function, where $\Tow(0,x)=x$ and
\begin{equation}\label{tower}
  \Tow(n,x)= e^{\Tow(n-1,x)} = \underbrace{e^{e^{\iddots^{e^x}}}}_{n
    \text{ iterations}}.
\end{equation}

\begin{thm}\label{T:k finite}
  Fix $\eps>0$, and consider the spatial Kingman coalescent on a graph $G$
  with bounded degrees. Start with $n$ particles located at $o\in G$, and
  let $N^n(t)$ be the total number of particles at time $t>0$. There are
  constants $C,c>0$ depending only on $t$ and the degree bound such that
  \begin{align*}
  \P\Big(N^n(t) \ge c \Vol B(o,(1-\eps)\log^* n) \Big)
  &\xrightarrow[n\to\infty]{} 1,  \qquad \text{ and } \\
  \P\Big(N^n(t) \le C \Vol B(o,(1+\eps)\log^* n) \Big)
  &\xrightarrow[n\to\infty]{} 1.
  \end{align*}
\end{thm}

A more concise way of stating \thmref{k finite} is that $N^n(t) \asymp \Vol
B(o,(1+o(1))\log^* n)$ with high probability (see the paragraph on ``Other notations''
at the end of Section \ref{S:defs} for $\asymp$ and other notations related to asymptotic behavior).

\begin{rem*}
  The function $\log^* n$ tends to infinity with $n$, but at a very slow
  rate: $\log^* n \le 4$ for $n \le 10^{1656520}$. Thus while this function
  diverges to infinity from the mathematically rigorous point of view, for
  all practical sample sizes it takes value $3$ or $4$.
\end{rem*}

\begin{rem*}
  The behavior in Theorem \ref{T:k finite} contrasts that of the mean-field
  case, where $N^n(t)$ converges (without renormalization) to a finite
  random variable $N(t)$ for all $t>0$, due to well-known properties of
  Kingman coalescent. In the lattice case $G=\Z^d$, we see that $N^n(t)$
  diverges as $(\log^* n)^d$, i.e.\ extremely slowly. Even on a regular
  tree, where balls have maximal volume given the degree, $N^n(t)$ diverges
  only as $e^{C\log^* n}$.
\end{rem*}

The mean-field $\Lambda$-coalescent processes
can be classified according to the
{\em coming down from infinity} (CDI) property.
For the partition-valued process $(\Pi_t,\, t\geq 0)$, this
means that the initial configuration $\Pi_0=\{\{i\}:\, i \in \N\}$
is countably infinite, but that
$\Pi_t$ contains only finitely many classes at any time $t>0$, almost surely.
The Kingman and the Beta-coalescents with parameter in a certain range
(see below) come down from infinity, while the
Bolthausen-Sznitman coalescent does not. If the spatial coalescent is
viewed as an interacting particle system, then CDI is
the property that the total number of particles in the system
at any fixed positive time remains bounded (tight), as the initial number
of particles tends to $\infty$.
It turns out that for the mean-field model (as well as the spatial model
with $G$ a finite graph),
either the total number of particles in the system
converges almost surely to a finite random variable, or it
diverges at any given time.

A natural question, and one of the motivations of this work, is whether a
similar dichotomy occurs in spatial coalescents. It is known that if
the mean-field $\Lambda$-coalescent comes down from
  infinity, then the number of particles in its spatial counterpart
will be locally finite (see
Section~\ref{S:defs}).
It is natural to ask whether the total number of particles can
nevertheless be infinite on infinite graphs.
We answer this question for general spatial $\Lambda$-coalescents.
It is remarkable that the answer is universal,
in that it does not depend on the driving measure $\La$ nor on
the geometry of the underlying infinite graph $G$.

\begin{thm}\label{T:global}
  For any measure $\Lambda$ on (0,1) and any infinite graph $G$, consider
  the spatial $\Lambda$-coalescent on $G$ started with $n$ particles at
  $o\in G$. If $N^n(t)$ denotes the total number of particles at time $t$, then
  $N^n(t)\to\infty$ almost surely, as $n\to\infty$.
\end{thm}

In particular, for $\Lambda$ such that the
mean-field coalescent comes down from infinity, the number of particles
will be locally finite, but globally infinite. For this reason we call this
phenomenon the \emph{global divergence} of the spatial $\Lambda$-coalescent.

Our next result on the fixed time asymptotics of $N^n$ concerns a
setting that is particularly relevant for some biological applications,
where the coalescence mechanism is given by the Beta-coalescent with
parameter $\alpha\in (1,2)$ (as defined in the next section).

\begin{thm}\label{T:beta}
  Fix $\alpha\in (1,2)$, and consider the spatial Beta($2-\alpha,\alpha)$
  coalescent on a graph $G$ with bounded degree. Start with $n$ particles
  located at $o\in G$, and let $N^n(t)$ be the total number of particles at time
  $t$. There are constants $C,c>0$ depending only on $t,\alpha$ and the
  degree bound such that
  \begin{align*}
    \P\Big(N^n(t) \ge c \Vol B(o,c\log\log n) \Big)
    &\xrightarrow[n\to\infty]{} 1,  \qquad \text{ and } \\
    \P\Big(N^n(t) \le C \Vol B(o,C\log\log n) \Big)
    &\xrightarrow[n\to\infty]{} 1.
  \end{align*}
\end{thm}

Similarly to \thmref{k finite}, this may be stated more concisely as $N_n(t) =
\Theta(\Vol B(o,\Theta(\log\log n)))$. When the graph has growth $\Vol
B(o,R)\asymp R^d$, this translates to $N^n(t) \asymp (\log\log n)^d$.

\medskip The above theorems describe the state of the system at a fixed
time $t$. We also provide estimates for the number of particles that
survive for a long time. Here the diffusion of particles plays a more
important role, hence the results depend in a more fundamental way on the
underlying graph. We focus on Euclidean lattices $G=\Z^d$, $d\geq 2$.

\begin{thm}\label{T:long_time}
  Assume that the coalescence mechanism is Kingman's coalescent. Let
  $G=\Z^d$, let $m=\log^* n$, and fix $\delta>0$. Then there exist some
  constants $c>0$ and $C>0$ (depending only on $d,\delta$) such that, if
  $d>2$,
  \[
  \P\left( cm^{d-2} < N^n(\delta m^2) < Cm^{d-2} \right)
  \xrightarrow[n\to\infty]{} 1,
  \]
  while, if $d=2$, then
  \[
  \P\left( c\log m < N^n(\delta m^2) < C\log m \right)
  \xrightarrow[n\to\infty]{} 1.
  \]
  If the coalescence mechanism is a Beta-coalescent with parameter
  $\alpha\in (1,2)$, then the same statement holds with $m= \log \log n$.
\end{thm}

One interpretation of this theorem is that, when the underlying graph $G$
is $\Z^d$, the resulting random particle system may also be thought of as a
microscopic description of the small-time evolution of a solution to the
parabolic nonlinear partial differential equation:
\begin{equation}\label{Smoluchowski}
\partial_t u = \frac12 \Delta u - \beta u^2,  \qquad  (\beta>0),
\end{equation}
starting from a singular initial condition, such as a Dirac delta measure
at a given spatial location. We refer the reader to \cite{reztwo,rezone} for a
discussion of this equation.

As suggested by \thmref{long_time},
for the study of the long-term particle system behavior
it is natural to rescale the particle
system's time by a factor of $m^2$, while rescaling space by a factor of $m$.
\thmref{long_time} indicated that the rescaled system should exhibit a
Boltzmann-Grad limiting behavior, i.e., the number of interactions
(intersections and coalescences) between a particle and all others over any
finite time interval is tight and does not tend to 0 as $n\to\infty$.
The behavior of the
rescaled system mirrors the system of Brownian coagulating particles
studied in  \cite{reztwo} and \cite{rezone}, in which the PDE
\eqref{Smoluchowski} is derived as the governing macroscopic behavior, and
is obtained as a particular case of the Smoluchowski system of PDEs. It is
worth pointing out that in the current case, the discrete structure of the
lattice remains important in
determining the frequency of coalescence events even after space and time
have been rescaled, which would alter the formula fixing the reaction
coefficient in the limiting PDE.

\begin{rem*}
  The $\log^*$ function featured in \thmref{k finite} might remind the
  reader of a result of Kesten \cite{KestenTPB, KestenJMB}, who
  studied the number of allelic types in a Wright-Fisher model with small
  mutation probability.

  In Kesten's model, allelic types take values in $\Z$, and the type of an
  offspring is identical to that of its parent, except on a mutation event
  of a small probability (inversely proportional to the total population
  size). When a mutation occurs, the offspring's type is chosen by adding
  an independent $\Z$-valued random variable (with some given, bounded,
  distribution) to the parent's type, that is, by making a random walk step
  from the parent's type. It turns out that the number of types (and, in
  fact, their relative positions in space) has an equilibrium distribution.
  It is shown in \cite{KestenTPB, KestenJMB} that the number of observed
  types at equilibrium is of order $\log^* n$, where $n$ is the sample
  size. The above Fisher-Wright model may seem closely related to the
  one-dimensional spatial Kingman coalescent, but on a closer look one
  realizes that the dynamics of the two models are quite different, and
  there is no direct relation between the results.

  Kesten's result may be phrased as follows. Let $T_n$ be the tree
  generated by Kingman's (non-spatial) coalescent started with $n$
  particles, and consider a branching random walk indexed by $T_n$. Then
  the number of distinct values at the leaves is of order $\log^* n$. A
  variation of the strategy used in Section~\ref{S:finitetimeKingman}
  applies in this setting, and can lead to an alternate proof of Kesten's
  result.
\end{rem*}

\subsection{Heuristics and proof ideas}
\label{S:ideas}
It is evident from Theorems~\ref{T:k finite} and \ref{T:beta} that the long
term behavior of the number of particles in the spatial coalescent depends
delicately on the precise nature of the coalescent. We now describe the
approximate behavior of the spatial coalescent started with a large number
of particles, all located at $o$. The proofs are mostly a detailed
treatment of the following heuristic observations.

To understand the finite initial condition, we turn to the infinite one.
Consider a given $\Lambda$-coalescent which comes down from infinity (see
below). Let $N_t$ be the number of particles in the (non-spatial)
coalescent started with $N_0=\infty$. For Kingman's coalescent it is the
case that $N_t\sim2/t$, whereas for Beta-coalescents with parameters
$(2-\alpha, \alpha)$ with $1<\alpha<2$, we have $N_t\sim c_\alpha
t^{1-\alpha}$ \cite{bbs2, bbs1}. The rough description that follows applies
to both of these, as well as more general coalescents. In general, one
would expect $N_t$ to be concentrated (for small $t$) around some function
$g(t)$ (such a function is found in \cite{bbl1}). The coalescent started
with $N$ particles is similar to the infinite coalescent observed from time
$g^{-1}(N)$ onward.

Consider now the non-spatial coalescent with emigration, where each
particle also disappears at some rate $\rho$. In fact, the parameter $\rho$
may depend on the size of the population, as long
as $n\rho(n)$ is non-decreasing. It turns out that for
coalescents that come down from infinity, the emigration does not influence
$N_t$ so much, and $N_t$ is still close to $g(t)$. The total number of
particles that emigrate when starting with $N$ particles is then close to a
Poisson variable with mean
\begin{equation}
  f(N) := \int_{g^{-1}(N)} g(t) \rho(g(t)) dt.
\end{equation}
(The upper bound of integration is some arbitrary constant.)

Now comes the key observation: if $N$ is large, the number of particles
migrating back into $o$ is negligible (under a technical condition that
holds for most spatial coalescents), and in fact, an overwhelming
proportion of those particles that emigrate will have emigrated by time
$g^{-1}(f(N))$. Thus we find that at this time, the number of particles at
$o$ and each of its neighbors is of order $f(N)$. A second observation is
that the resulting populations can be approximated by independent spatial
coalescents, when observed from time $g^{-1}(f(N))$ onward. In particular,
at time $g^{-1}(f\circ f(N))$ there are of the order of $f\circ f(N)$
particles at each vertex in $B(o,2)$. This ``cascading onto neighbors''
continues until step $m$, where $m$ is such that $f \circ \dots \circ f(N)$
($m$ repeated iterations of $f$) is of order 1. Note that in these $m$
steps a ball of radius $m$ has been roughly filled.

Applying this heuristics to the case of Kingman's coalescent and the
Beta-coalescents with parameters $(2-\alpha, \alpha)$ and $1<\alpha <2$,
gives the following. For Kingman's coalescent and constant $\rho$, we have
$f(n)\sim 2\rho\log n$, and for Beta-coalescents we have $f(n)\sim C_\alpha
\rho n^{2-\alpha}$, for some constant $C_\alpha>0$. Thus in the first case,
$m=\log^* N$. In the second case, we find $m\sim c\log\log N$. In general,
this gives $m \sim f^*(n)$, where
\begin{equation}\label{fstar}
  f^*(n) = \inf \bigg\{ m \ge 1: \underbrace{f \circ \dots \circ f}_{m\text{
        iterations}}(n) \le 1 \bigg\}.
\end{equation}
Note that if $\rho(n)$ decreases fast enough so that $f(n)$ is bounded,
then it follows from this heuristic analysis that the spatial coalescent
will come down from infinity globally. However, when $\rho$ is constant, it
can be proved that $f$ is always unbounded, which in turn implies the
result about global divergence of any spatial $\Lambda$-coalescent.

\medskip

Turning to the long time asymptotics, by the above reasoning we may start from
a configuration consisting of a tight number of particles at each site of
the ball of radius $m$ around the origin. Since the number of particles per
site is tight, the coalescent dynamics influences the evolution less than
the diffusion. In particular, the structure of the underlying graph becomes important for
the asymptotic behavior of the process. For simplicity, let us restrict
ourselves to $d$-dimensional Euclidean lattices with $d\ge 3$. Let
$\rho(t)$ denote the average number of particles per site in the ball of
radius $m$ at time $t$. Then at time $t_0=1$ we have $\rho(t_0)\asymp 1$
and $\lim_{t\to \infty} \rho(t)=0$. Each particle present in the
configuration at time $t$ coalesces with another particle at an average
rate approximately $\rho(t)$, so that $\frac{d}{dt} N(t) = - N(t)
\rho(t)/2$. Dividing by the volume of the ball, one arrives to the ODE
\begin{equation}\label{rho_heur}
  \frac{d}{dt} \rho(t) = - \frac12\rho(t)^2 ,
\end{equation}
whose solution is given by $\rho(t)=2/(t+c)$ for some $c>0$.

The approximation (\ref{rho_heur}) should be valid as long as the diffusion
of particles away from the initial region (i.e., $B(o,m)$) is negligible.
The influence of diffusion should start to be visible at times of order
$m^2$. In particular, at time $m^2$, the density $\rho(m^2)$ is of order
$m^{-2}$, so the total number of remaining particles is of order $m^{d-2}$.
Assuming the plausible claim that the remaining particles are approximately
uniformly distributed over a ball of radius order $m$, a simple calculation
(using hitting probabilities for random walks) now implies that each of
them has a positive probability of never meeting any other particle again,
and so the number of particles that survive indefinitely is of order
$m^{d-2}$.

We wish to point out that van den Berg and Kesten \cite{BKcoal1, BKcoal2} have
shown a density decay similar to (\ref{rho_heur}) for a related model of
coalescing random walks. However their results differ in two ways. On the
one hand, the coalescence mechanism which they analyze is different. On the
other hand, and more importantly, their initial condition is initially
homogeneous in space, and not restricted to a large ball. This restriction
is the cause of much of the difficulty in the current setting -- see Section
\ref{SS:upper bound} for more details.

\subsection{Definitions and background on spatial coalescents}
\label{S:defs}

\paragraph{Kingman's coalescent.}
Suppose that we are given an integer $n\ge 1$. {\em Kingman's $n$-coalescent} is
the Markov process $(\Pi^n_t,t\ge0)$, with values in the set
$\mathcal{P}_n$ of partitions of $[n]:=\{1,\ldots,n\}$, such that $\Pi^n_0=
\{\{1\},\{2\},\ldots,\{n\}\}$, and such that each pair of blocks merges at
rate 1, and these are the only transitions of the process. Blocks of the
partition $\Pi^n_t$ may be viewed as indistinguishable particles, and we
often refer to the number of blocks of $\Pi^n_t$ as the number of particles
alive at time $t$. A simple but essential property of Kingman's
$n$-coalescent is the so-called {\em sampling consistency} property: the
restriction of $(\Pi^{n+1}_t,t\ge 0)$ to $[n]$ has the same distribution as
an $n$-coalescent. This enables one to construct a Markov process
$(\Pi_t,t\ge 0)$ with state space $\mathcal{P}$, the set of partitions of
$\N$, such that the law of $\Pi$ when restricted to $[n]$ equals the law of
$\Pi^n$. In particular, the initial state of this process is the trivial
partition $\Pi_0=\{\{1\}, \{2\},\ldots\}$. The process $\Pi$ is called {\em
  Kingman's coalescent}. For background reading, see for instance
\cite{durrettbookDNA, saintflour, coalnotes}.

\paragraph{$\Lambda$-coalescents.}
Let $\Lambda$ be a finite measure on $[0,1]$. A {\em coalescent with
  multiple collisions}, or {\em $\Lambda$-coalescent}, is a Markov process
$(\Pi_t,t\ge 0)$ with values in the set of partitions of $\N$ characterized
by the following properties. If $n \in \N$, then the restriction of
$(\Pi_t,t\ge 0)$ to $[n]$ is a Markov chain $(\Pi^{(n)}_t,t\ge 0)$, where
$\Pi^n_0= \{\{1\},\{2\},\ldots,\{n\}\}$, and where the only possible
transitions are mergers of blocks (it is possible to merge several blocks
simultaneously into one block, but no two mergers of this kind can occur
simultaneously) so that whenever the current configuration consists of $b$
blocks, any given $k$-tuple of blocks merges at rate
\begin{equation}\label{rate coal}
  \lambda_{b,k}=\int_{[0,1]} x^{k-2}(1-x)^{b-k}\Lambda(dx).
\end{equation}
Note that $0^0$ is interpreted as 1, so that an atom of $\Lambda$ at 0
causes each pair of particles to coalesce at a finite positive rate
$\La(\{0\})$. In this way any $\La$-coalescent can be thought of as a
superposition of a ``pure'' coalescent with multiple collisions driven by
measure $\La(dx)\indica(0,1](x)$, and a time-changed Kingman's coalescent.
An atom of $\La$ at 1 causes all the particles to coalesce at some positive
fixed rate. Such $\La$-coalescent may be viewed as a killed
$\La'$-coalescent where $\La'(dx)=\La(dx)\indica[0,1)(x)$. Kingman's
coalescent is a particular $\Lambda$-coalescent, obtained when the measure
$\Lambda$ equals $\delta_{0}$, the unit Dirac mass at $0$. Any
$\Lambda$-coalescent $\Pi$ is sampling consistent, that is, if $m<n$ then
the restriction of $\Pi^n$ to $[m]$ is equal in law to $\Pi^m$. It is this
observation that allows one to construct an infinite version of the
process. It is interesting to note the following fact shown by Pitman
\cite{pit99}: $\Lambda$-coalescents are the only exchangeable Markov
coalescent processes without simultaneous collisions. We refer the reader
to \cite{pit99} for definitions and further properties.

As already mentioned,
if $\Lambda(dx)=dx\indica{[0,1]}$, the corresponding $\La$-coalescent is
usually called the Bolthausen-Sznitman coalescent, and more generally if
$\Lambda$ is the
Beta($2-\alpha,\alpha$) distribution where $\alpha\in (0,2)$ is a fixed parameter, that is,
\begin{equation} \label{beta}
  \Lambda(dx) = \frac{1}{\Gamma(2-\alpha)
    \Gamma(\alpha)} x^{1-\alpha} (1-x)^{\alpha - 1} \: dx,
\end{equation}
the corresponding $\Lambda$-coalescents is called {\em Beta-coalescents} with parameter $\alpha$.
The Bolthausen-Sznitman coalescent is the special case $\alpha=1$, and it does not come down from infinity.
For $\alpha \in (1,2)$, the corresponding Beta-coalescents come down from infinity,
and they are important processes from the theoretical evolutionary biology perspective,
due to the following result from \cite{schweinsberg}: the Beta-coalescent with parameter
$\alpha\in (1,2)$ arises in the scaling limit of population models where the
offspring distribution of a typical individual is in the domain of
attraction of a stable law with index $\alpha$. Apart from the Kingman
coalescent, the Beta-coalescents with parameter $\alpha\in (1,2)$
are the most-studied class of $\La$-coalescents (see, e.g., \cite{birkner-etal, bbs1, bbs2}).

\paragraph{Spatial coalescents.}
As informally described above, spatial coalescents are processes which
combine spatial motion of individual particles with coalescence of
particles located on the same site of a given graph of bounded degree.
Let $\Lambda$ be a given finite measure on $[0,1]$. A spatial
$\Lambda$-coalescent, as defined in \cite{LS}, is a Markov processes
$(\Pi^\ell_t,t\ge 0)$ with values in the space $\mathcal{P}^\ell =
\mathcal{P}\times V^{\{1,2,\ldots\}}$ of partitions of $\{1,2,\ldots\}$
indexed by spatial locations. That is, an element $x=(\pi,\ell)\in
\mathcal{P}^\ell$ consists of a partition $\pi=\{A_1,A_2,\ldots\}$, and a
sequence $\ell=(\ell_1,\ell_2,\ldots)$, where $\ell_i$ specifies the
location of the block $A_i$. There are only two types of transitions
possible for $\Pi^\ell_t=(\Pi_t,\ell_t)$: (i) provided there are $b$ blocks
at a location $v\in V$, then any given $k$-tuple of them will merge at rate
$\lambda_{b,k}$ given by \eqref{rate coal}, independently over $v$; and
(ii) independently of the coalescent mechanism, each block $A_k$ of $\pi$
migrates at rate $\theta$. This means that if the block is at $v$, then
some vertex $w$ is chosen according to the distribution $p(v, \cdot)$,
where $p(v,w)$ is a given Markov kernel. When this happens, $\ell_k$ is
changed from $v$ to $w$. To simplify the discussion, we will assume unless
otherwise specified, that $p(x,y)$ is the transition kernel for the simple
random walk on the underlying graph $G$.

If $\pi$ is a partition let $i\sim_\pi j$ mean that the particles labeled
$i$ and $j$ belong to the same block of $\pi$. For $(\pi,\ell) \in
\mathcal{P}^\ell$ and $v \in V$, denote by $\#_v(\pi,\ell)$ the number of
blocks in $\pi$ with label (location) $v$.

Spatial $\La$-coalescents inherit the sampling consistency directly from
$\Lambda$-coalescents. Namely, if we consider a spatial coalescent started
from $n+1$ particles (that is, blocks) and consider its restriction to the
first $n$ particles, the new process has the law of a spatial coalescent
started from $n$ particles. This simple property will be used on several
occasions. In particular, it implies that if $(\pi^1,\ell^1)$ and
$(\pi^2,\ell^2)$ are such that $\#_v(\pi^1,\ell^1)\leq \#_v(\pi^2,\ell^2)$,
for all $v$, then there exists a coupling of two spatial coalescents
$((\Pi_t^1,\ell_t^1),(\Pi_t^2,\ell_t^2)), t\geq 0)$ such that
$(\Pi^i_0,\ell_0^i)=(\pi^i,\ell^i)$, $i=1,2$ and $\#_v
(\Pi_t^1,\ell_t^1)\leq \#_v (\Pi_t^2,\ell_t^2)$ for all $v$, almost surely.
The same property guarantees the existence of spatial coalescents started
with infinitely many particles on an infinite graph (see Theorem 1 in
\cite{LS} for a particular construction).

Spatial $\La$-coalescents may be started from configurations containing
countably infinitely many particles at each site of $G$, see \cite{LS}.
However, our main results concern spatial $\La$-coalescents started from the
following initial condition:
\begin{equation}\label{IC}
  \Pi^\ell_0=(\{\{1\},\{2\},\ldots\}, (o,o,\ldots)),
\end{equation}
where $o$ is some given reference vertex called the {\em origin} of $G$. In
words, all the infinitely many particles are initially located at the
origin $o$.

From now on we abbreviate
\begin{equation}
  \label{eq:nota Xv}
  X_v(t) = \#_v (\Pi_t,\ell_t) \qquad \text{ and } \qquad
  X_v^n(t) = \#_v (\Pi^n_t,\ell_t).
\end{equation}
We denote the {\em total number of blocks} by $\tN(t)=\sum_{v\in V} X_v(t)$
(resp.~$N^n(t)=\sum_{v\in V} X^n_v(t)$). When not in risk of confusion, we
will drop the superscript $n$ to simplify notations. It is clear from the
definitions that both $(\{X_v(t)\}_{v\in G},t\geq 0)$ and $(\tN(t),t\geq
0)$ have Markovian transitions, with respect to the filtration generated by
the coalescent process $\Pi$. They carry only partial information about the
evolution of the corresponding spatial coalescent, in particular, they do
not determine the evolving partition structure.

In the language of theoretical population biology, a sample of $n$
individuals is selected from the population at the present time, and their
ancestral lineages are followed in reversed time. The above
transition rules (i)--(ii) given above reflect the idea that individuals typically reproduce
within their own colony (so that only particles on the same site may
coalesce), and occasionally there is a rare migration event, which
corresponds to the random walk transitions. In the case where the
coalescence mechanism is simply Kingman's coalescent, we note that this
model may be viewed as the ancestral partition process associated with
Kimura's {\em stepping-stone model} \cite{Kimura, KimuraWeiss}.

\paragraph{Coming down from infinity.}
Let $(\Pi_t,t\ge 0)$ be Kingman's coalescent. As already mentioned, Kingman
\cite{king82, king82b} realized that while $\Pi$ starts with an infinite
number of blocks at $t=0$, its number of blocks becomes finite for all
$t>0$, almost surely. A coalescent with multiple collisions
 may or may have the same property, depending on the measure $\Lambda$. More
precisely, there are only two possibilities as shown in \cite{pit99}:
let $E$ (resp.~$F$) denote the event that for all $t>0$
there are infinitely (resp.~finitely) many blocks. Then, if
$\Lambda(\{1\})=0$, either $P(E)=1$ or $P(F)=1$. When $P(F) = 1$, the
process $\Pi$ is said to {\em come down from infinity}. For instance, a
Beta-coalescent comes down from infinity if and only if $1<\alpha < 2$,
henceforth we make this an assumption whenever working with
Beta-coalescents.

In the context of spatial coalescents, assuming that $\Lambda(\{1\})=0$,
Proposition 11 in \cite{LS} implies that when the initial number of
particles is infinite, then $X_v(t)$ becomes finite for all $v\in V$ and
$t>0$ with probability 1, if and only if the underlying measure $\Lambda$
is such that the mean-field (i.e., non-spatial) $\Lambda$-coalescent comes
down from infinity. In this situation, we may say that the spatial
coalescent comes down from infinity {\em locally}. Naturally, this stays
true if the initial condition is \eqref{IC}.

\paragraph{Other notations.}
Unless specified otherwise, $c,C$ (and variations $c_1,C_2,\ldots$) will
henceforth denote positive constants that depend only on the underlying
graph,
and that may change from line to line. Typically, $c,c_1,\ldots$ denote
sufficiently small, whereas $C,C_1,\ldots$ denote sufficiently large
constants. We also use the symbols $a_n\sim b_n$ and $a_n\asymp b_n$ to
denote respectively that $a_n/b_n \to 1$, and $a_n/b_n$
is bounded away from 0 and $\infty$, as $n \to \infty$.

\paragraph{Organization of the paper.}
The rest of the paper is organized as follows. Section~\ref{S:prelim}
starts with some preliminary remarks and observations concerning
large deviation estimates for Kingman's
coalescent and Ewens's sampling formula, as well as
several couplings between the spatial $\Lambda$-coalescents and the corresponding
(mean-field) $\Lambda$-coalescents, which will be used throughout
the paper. Section~\ref{S:finitetimeKingman} contains a proof of \thmref{k
  finite} on the behavior of the spatial Kingman coalescent in finite
time. As many of the subsequent results in the paper build on this, we
recommend reading this section prior to any of the following sections.
Section~\ref{S:ProofBeta} contains
the proof of \thmref{beta} on the finite-time behavior of the spatial Beta-coalescents.
Section~\ref{S:global} returns to the general
case of spatial $\Lambda$-coalescents and arbitrary graphs with bounded
degree, and contains the proof of global divergence (\thmref{global}). In
final Sections~\ref{sec:long time} and \ref{SS:upper bound} we study
respectively the lower bound and the upper bound for long term
behavior of Kingman's coalescent (as stated in \thmref{long_time}). The lower bound
 obtained in Section~\ref{SS:lowerbound2} is true for general $\Lambda$-coalescents,
but we provide in Section~\ref{SS:concentration} an alternate shorter proof
for the special case of Kingman's coalescent, that also gives tighter
bounds. The proof of the upper bound in Section~\ref{SS:upper bound} turns
out to be the most technical part of the paper, and it is based on a delicate multi-scale
analysis.

Sections \ref{S:ProofBeta}--\ref{sec:long time} may be read in any order, depending on
the interest of the reader. We recommend reading Section \ref{sec:long time} prior to
Section \ref{SS:upper bound}.

\section{Preliminary lemmas}
\label{S:prelim}

\subsection{Some large deviation estimates}

We begin with an easy Chernoff type bound for a sum of exponential random
variables, which we prefer to state in an abstract form now so as to refer
to it on several occasions later. In our applications, $\E S$ will
typically be small.

\begin{lemma}\label{L:ldev}
  Let $\{E_i\}_{i\in I}$ be independent exponential random variables with
  $\E E_i = \mu_i$. Let $S = \sum_{i\in I} E_i$. Then for any $0<\eps<1$
  \[
  \P\big( S < (1-\eps) \E S \big) \le \exp\left(-\frac{\eps^2(\E
      S)^2}{4\Var S} \right).
  \]
  Additionally, for $0 < \eps < \frac{\Var S}{\E S \sup\{\mu_i\}}$,
  \[
  \P\big( S > (1+\eps) \E S \big) \le \exp\left(
    - \frac{\eps^2 (\E S)^2}{4\Var S}\right).
  \]
\end{lemma}

\begin{rem*}
  If $I=\{n,n+1,\dots\}$ and $\mu_i \sim c i^{-\alpha}$ for some
  $\alpha>1$, then as $n\to \infty$, $\frac{\Var S}{\E S \sup\{\mu_i\}}$ is
  bounded away from 0, hence the second bound holds for all $\eps>0$ small
  enough, for all $n$.
\end{rem*}

\begin{proof}
  Using Markov's inequality, for any $0 < \lambda \le \frac12 \inf
  \{\mu_i^{-1}\}$
  \begin{align*}
    \P\big( S > (1+\eps) \E S \big)
    &\le e^{-\lambda(1+\eps)\E S} \E e^{\lambda S} \\
    &= e^{-\lambda(1+\eps)\E S} \prod \frac1{1-\lambda \mu_i} \\
    &< e^{-\lambda(1+\eps)\E S}
    \exp\left(\sum \lambda\mu_i + \lambda^2\mu_i^2 \right) \\
    &= e^{- \lambda \eps \E S + \lambda^2 \Var S},
  \end{align*}
  where we have used that for $x\in(0,1/2)$ we have $-\ln(1-x) < x+x^2$.
  Taking $\lambda = \frac{\eps\E S}{2\Var S}$, which is allowed since $\eps
  < \frac{\Var S}{\E S \sup\{\mu_i\}}$, yields the upper bound.

  The lower bound follows from a similar argument with $\lambda =
  -\frac{\eps\E S}{2\Var S}$.
\end{proof}

We now apply this to get a large deviation estimate for Kingman's
coalescent. This uses a simple idea which can already be found in Aldous
\cite{aldous}, who used it to prove a central limit theorem for the number
of particles at time $t$. Denote by $\P^n$ the law of the (non-spatial)
Kingman coalescent started with $n$ blocks. Let $N(t)$ be the number of
blocks at time $t$.

\begin{lemma}\label{L:kingman ldev}
  Let $t=t(n)\to 0$ in such a way that $t(n)^{-1}=o(n)$. For any $0 < \eps
  < 1/2$, for $n$ large enough,
  \[
  \P^n\left( 1-\eps < \frac{N(t)}{2/t} < 1+\eps \right) > 1 - \exp\left(
    -\frac{\eps^2}{t} \right).
  \]
\end{lemma}

\begin{proof}
  For the upper bound, let $m=\lceil (1+\eps)2/t \rceil$. The time it takes
  the process to get from $n$ to $m$ particles is a sum of independent
  exponential random variables with means $\binom{k}{2}^{-1}$ for
  $k=m+1,\dots,n$. Call this sum $S$. If $N(t)>m$ then $S>t$. We have
  \[
  \E S = \sum_{k=m+1}^n \binom{k}{2}^{-1} \sim 2m^{-1} \sim t/(1+\eps)
  \]
  provided $m=o(n)$. Similarly,
  \[
  \Var S = \sum_{k=m+1}^n \binom{k}{2}^{-2} \sim (4/3) m^{-3}.
  \]
  Thus, for $\eps<2/3+o(1)$,
  \[
  \P(S > t) < \exp\left(-\frac{(3+o(1))\eps^2}{2t} \right),
  \]
  by \lemref{ldev}. The lower bound is similar using the upper bound on
  $S$.
\end{proof}

We now consider Kingman's coalescent with spatial migration. Let $\P^n$ be
the law of a simplified process where $n$ particles initially located at a
single site $o$ coalesce according to Kingman's dynamics, while each
particle (or block of particles) migrates at rate $\rho$, and any block
that migrates away from $o$ is ignored from that time onwards. Denote by
$Z_n$ the total number of blocks that ever migrate away from $o$.

\begin{figure}
  \begin{center}
    \includegraphics[width=6cm]{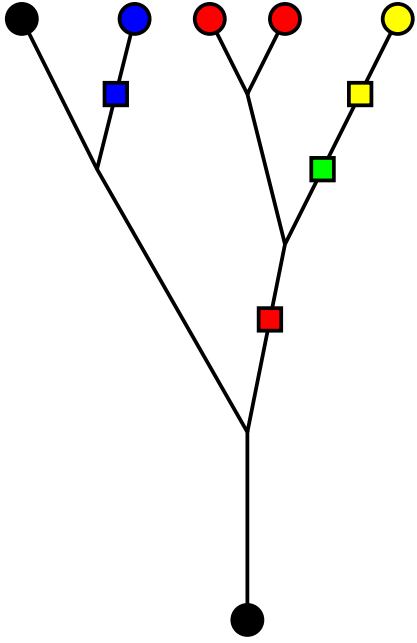}
    \caption{The random partition generated by mutations (squares). Here
      $Z_5=4$. }
    \label{Fig:ESF}
  \end{center}
\end{figure}

One can think of each migration event as of a ``unique mutation on the
genealogical tree'', by giving it for example the label equal to its occurrence time.
Since migrations happen at rate $\rho$ for each block present in the
configuration at site $o$, one quickly realizes that $Z_n$ is a realization
from a well-known distribution arising in mathematical population genetics.
Namely, set $\theta=2\rho$, and suppose that on the (non-spatial) Kingman
coalescent tree mutation marks occur at a Poisson rate of $\theta/2$ per
unit length. Using the language of mathematical population genetics, assume
the infinite alleles models (all mutations create a different allele, and
so different individuals in the original sample of $n$ are in the same
family if and only if they descend from the same mutation and there has
been no other mutation between this common ancestor and the present
individuals). The marks of the mutation process generate a random
partition $\Pi_\theta$ on the leaves of the tree by declaring that $i$ and
$j$ are in the same block of $\Pi_\theta$ if and only if there is no
mutation mark on the shortest path that connects $i$ and $j$. In
Figure~\ref{Fig:ESF} different blocks of this partition are represented by
different colors. Then it is easy to see that $Z_n$ has the law of the
number of blocks in $\Pi_\theta$. It is well-known (see, e.g., (3.24) in
Pitman \cite{saintflour}) that $Z_n$ is of order $\theta \log n$ for large
$n$. The following large deviation estimate is part of the folklore, but we
could not find a precise reference for it in the literature.

\begin{lemma}\label{L:esf ldev}
  Fix $\eps>0$. There are $c,C>0$ such that
  \[
  \P^n \left(\left|\frac{Z_n}{\log n} - \theta \right| > \eps \right) <
  Cn^{-c}.
  \]
  Furthermore, for any $U$,
  \[
  \P^n(Z_n > U) < C n^C e^{-U}.
  \]
\end{lemma}

The proof is based on the Chinese restaurant process representation of
Ewens's sampling formula. Let $K_{n,i}$ be the number of blocks of size $i$
in $\Pi_\theta$, where $i=1,\ldots,n$. Then the distribution of
$(K_{n,1},\ldots,K_{n,n})$ is given by Ewens's sampling formula
$ESF(\theta)$:
\begin{equation}
  \label{ESF}
  P(K_{n,i}=a_i,\,i=1,\ldots,n) =
  \frac{n!}{\theta (\theta+1)\cdots(\theta+n-1}\prod_{i=1}^n
  \frac{\theta^{a_i}}{i^{a_i}a_i!}\,,
\end{equation}
for any given collection $a_1,\dots,a_n$ of non-negative integers such that
$\sum_{i=1}^n ia_i=n$.

The Chinese restaurant process representation of (\ref{ESF}) (see \cite[\S
3.1]{saintflour}), states that the number of blocks in $\Pi_\theta$
satisfies
\begin{equation}\label{eq:zetasum}
  Z_n \overset{d}{=} \sum_{i=1}^n \zeta_i
\end{equation}
where $\zeta_i$ are independent Bernoulli random variables with mean
\[
P(\zeta_i=1)= \frac{\theta}{i+\theta}.
\]

\begin{proof}[Proof of \lemref{esf ldev}]
  By \eqref{eq:zetasum} we have for $\lambda>0$
  \begin{align*}
    \E e^{-\lambda Z_n} &= \prod_{i\leq n}
    \left( 1 - \frac{\theta}{i+\theta}(1-e^{-\lambda})\right) \\
    &< \exp\left( \sum_{i\leq n}
      - \frac{\theta}{i+\theta}(1-e^{-\lambda}) \right) \\
    &< \exp\left( -\theta(1-e^{-\lambda})(C+\log n) \right).
  \end{align*}
  By Markov's inequality
  \begin{align*}
    \P(Z_n < (1-\eps)\theta\log n) &< \exp\left(\lambda(1-\eps)\theta\log n
      - \theta(1-e^{-\lambda})(C+\log n) \right) \\
    &< \exp\left((-\lambda \eps \theta + O(\lambda^2))\log n + C\right)
  \end{align*}
  For small positive $\lambda$ the coefficient of $\log n$ is strictly
  negative.

  Similarly, for $\lambda>0$
  \begin{align*}
    \E e^{\lambda Z_n}
    &= \prod_{i=1}^n \left(1+\frac{\theta}{i+\theta}(e^\lambda-1)\right) \\
    &< \exp\left( \sum_{i=1}^n \frac{\theta}{i}(e^\lambda-1) \right) \\
    &< \exp\Big( \theta(e^\lambda-1)(C+\log n) \Big).
  \end{align*}
  Thus, by Markov's inequality,
  \begin{equation}\label{eq:Zn not large}
    \P \left( Z_n > U \right)
    < \exp\Big( \theta(e^\lambda-1) (C + \log n) -\lambda U \Big)
  \end{equation}
  Taking $U=(1+\eps)\theta\log n$ and $\lambda=\lambda(\eps,\theta)$ small
  enough gives the first upper bound. Taking $\lambda=1$ gives the second
  claim.
\end{proof}

A similar computation can be found in Greven et al.\ \cite[Lemma~3.3]{GLW}.

\subsection{Coupling and comparison}

The final general tool we use is a coupling of the spatial process with
simpler coalescent processes such as non-spatial ones. While it is usual in
coalescent theory to keep track of the entire partition structure as time
evolves, here we are only interested in the number of particles remaining
in the system at any given time, so accounting for the full partition data
is cumbersome.

However, the detailed {\em Poisson process construction} (e.g., \cite{LS} Theorem 1)
becomes useful in the context of coupling.
Its advantage comes from the following fact:
all the information on both the merging and the migration is
given by the Poisson clocks
(ringing for jumps and for mergers), hence
one builds the spatial
coalescent process path (keeping track of the particle labels, and of
the partition structure) by applying a deterministic
function $\omega$-by-$\omega$ to this data.
Therefore, it is straightforward to append
another deterministic ingredient (as the ``coloring procedure'' in the following
lemma) to the construction.

Henceforth, it is convenient to consider the following simpler variation,
where the partition structure is ignored. Label the $n$ initial particles
by $1, \ldots, n$, and let $x_1, \ldots, x_n$ be their initial locations.
Let $S^1, \ldots, S^n$ be $n$ i.i.d.~simple random walks on $G$ in
continuous time with jump rate $\rho$ started at $x_1, \ldots, x_n$,
respectively. To each $k$-tuple of labels $i_1< \ldots< i_k$, and each
$b\geq k$, corresponds an independent Poisson process
$M_{i_1,\ldots,i_k}^b$ with intensity $\lambda_{b,k}$. The particles
labeled $i_1<\dots<i_k$ coalesce at a jump time $t$ of
$M^b_{i_i,\dots,i_k}$ if and only if they are all located at the same site
$v$ at time $t^-$, and there are a total of $b$ particles at $v$. In this
case, the newly created particle inherits the minimal label $i_1$.
(Subsequently there are no particles with labels $i_2,\dots,i_k$.) In
particular, its trajectory starting from time $t$ will be
$(S_s^{i_1},\,s\geq t)$.

Suppose that the $n$ initial particles are partitioned into classes
according to a partition $\pi=(B_1, \ldots, B_r)$ of $\{1, \ldots, n\}$, where
$r\ge 1$. We wish to compare the system $X=\{X_v(t),t\ge 0\}_{v \in V}$ to
the one which consists only of the particles that belong to a particular class
$B$ of $\pi$.
More precisely, for each $B\in \pi$, denote by $X^B$ the spatial coalescent process
whose initial configuration contains only the particles from $B$.
Denote by $N(t)$ the total number of particles of $X(t)$ and by
$N^B(t)$ the total number of particles of $X^B(t)$.

\begin{lemma}\label{L:sum_dom}
  There is a coupling $(X,X^{B_1}, \ldots, X^{B_r})$ such that, almost
  surely,
  \[
  \forall v\in V \ \ X^{B}_v(t) \leq X_v(t) \leq \sum_{i=1}^r X^{B_i}_v(t),
  \]
  for each block $B$ of $\pi$, and hence
  \[
  \ N^{B}(t) \leq N(t) \leq \sum_{i=1}^r N^{B_i}(t).
  \]
\end{lemma}

Note that in the coupling given below, the processes $X^{B_i}$ are not
independent. In fact, under weak assumptions on the coalescent and when all
the blocks are ``small'', there is a coupling including independence, see
\lemref{sumdom2}.

\begin{proof}
  Fix a realization of the process $(X(t),t\ge 0)$ as described above. Note
  that at any time $t$, each particle in the current configuration can be
  identified with a set of particles from the original configuration, that
  have been merging (possibly in several steps) to form this particle (this
  is the partition-valued realization of the spatial coalescent). In the
  rest of the argument, we say that a particle intersects $B \in \pi$ if
  its corresponding set intersects $B$.

  For each $i=1,\ldots,r$, and $t\geq 0$, let $X^{B_i}(t)$ be the
  configuration obtained from $X(t)$ by restricting to only those particles
  that intersect $B_i$. The consistency property of $\Lambda$-coalescents
  implies that for each $i$ the law of $X^{B_i}$ is that of the spatial
  coalescent started from the initial configuration restricted to elements
  of $B_i$. Thus this is a coupling of the processes.

  In this construction we have $X^{B_i}_v(t) \le X_v(t)$, since $X_v$ may
  contain particles that do not intersect $B_i$. Moreover, any particle
  contributing to $X_v(t)$ intersects at least one $B_i$, giving the bound
  $X_v(t) \leq \sum_i X^{B_i}_v(t)$. The inequalities relating $N(t)$ and
  $N^{B_i}(t)$ are an immediate consequence.
\end{proof}

A second type of coupling we will need is between a spatial
$\Lambda$-coalescent and its mean-field (i.e., non-spatial) counterpart. Fix a vertex
$u \in V$ of the graph, and consider a spatial $\Lambda$-coalescent
$\{X_v(t),t\ge 0\}_{v \in V}$ started with a finite number of particles and
such that initially $X_u(0) = n.$ Let $M(t)$ denote the number of particles
on $u$ at time $t$ that have always stayed at $u$, and let $Z(t)$ denote
the number of particles that jumped out of $u$ prior to time $t$. In
parallel, let $(N(t),t\ge 0)$ denote the number of particles at time $t$ in
a mean-field $\Lambda$-coalescent started with $n$ particles.

\begin{lemma}\label{L:coupl}
  There exists a coupling of $X$ and $N$ such that:
  \begin{equation}\label{EMZN}
    M(t) \leq N(t) \leq M(t)+Z(t), \qquad \text{ a.s.\ for all $t\geq 0$,}
  \end{equation}
  and
  \begin{equation}\label{eq:X_NZ}
    N(t)-Z(t) \leq X_u(t) \leq N(t) + Z(t), \mbox{ a.s.\ for all $t\geq 0$}.
  \end{equation}
\end{lemma}

\begin{proof}
  The process $M(t)$ may be realized as a mean-field coalescent where, in
  addition, particles are killed at rate $\rho$. In that case, if we let $
  Z(t)$ denote the total number of particles that have been killed, we see
  immediately that on the one hand, $M(t) \le N(t)$, and on the other hand,
  $N(t) = M(t) + \bar Z(t)$ where $\bar Z(t) \le Z(t)$. Indeed, $M(t)+
  Z(t)$ counts the number of particles if we freeze particle instead of
  killing them. However, in $N(t)$ these particles keep coalescing, and so
  the difference $\bar Z(t) = N(t) - M(t) \le Z(t)$. This proves
  (\ref{EMZN}). For (\ref{eq:X_NZ}), note first that $X_u(t) \ge M(t) =
  N(t) - \bar Z(t) \ge N(t) - Z(t)$. Finally, the last inequality in
  (\ref{eq:X_NZ}) is obtained by observing that $X_u(t)$ is made of
  particles that never jumped out of $u$ (there are $M(t)$ such particles)
  and of particles that have jumped out of $u$ and have come back at some
  time later, potentially coalescing in the meantime. There can never be
  more than $Z(t)$ such particles, since this is the total number of
  particles that jump out of $u$.
\end{proof}

In fact, one can be slightly more precise than the above estimate. We shall
need the following observation. Define two processes
\begin{align}\label{E mart S}
  S(t) &= Z(t) - \int_0^t \rho M(s) ds & V(t) &= S(t)^2 - \int_0^t \rho
  M(s) ds.
\end{align}
It is a standard (and easy) fact that both are continuous time martingales
under the law $\P^n$, with respect to the filtration $\FF$ generated by the
above coupling process. In fact, if we define $\GG=\sigma\{N(u), u\geq 0\}$
to be the $\sigma$-algebra generated by $N$, and let $\FF^*_t=\sigma
\{\GG,\FF_t\}$, then the processes $S(\cdot)$ and $V(\cdot)$ are
continuous-time martingales with respect to the filtration $\FF^*$.

\begin{lemma}\label{L:few_leave}
  For each time interval $[a,b]$, we have the stochastic domination
  \[
  \P(Z(b)-Z(a)\geq x|\GG) \leq \P\left(\Poi \left( \left.\rho \int_a^b N(s)
        \,ds \right)\geq x\right|\GG\right).
  \]
\end{lemma}

\begin{proof}
  Given $\GG$, $Z$ is a pure jumps process with jumps of size $1$ that
  arrive at rate $\rho M(t)\leq \rho N(t)$ at time $t$, almost surely.
 \end{proof}

 Finally, a global comparison with mean-field coalescents can be obtained
 in the case of the spatial Kingman coalescent as follows (see also
 \cite[\S6.1]{GLW1}). Let $S$ be an arbitrary subset of vertices and
 consider the restriction of $X$ to $S$.

\begin{lemma}\label{L:k-few in S}
  Fix a time $\tau\leq2$, and vertex set $S$, and assume that all particles
  are in $S$ at time 0. Let $Z=Z(\tau)$ be the number of distinct particles
  that exit $S$ by time $\tau$, and let $N_S(t)$ be the number of particles
  in $S$ at time $t$. Then, for some $c = c(\epsilon) >0$,
  \[
  \P\left( N_S(\tau) > Z + \frac{(4+\eps)|S|}{\tau} \right) <
  e^{-c|S|/\tau}.
  \]
\end{lemma}

Note that the bound is independent of the starting configuration. This
lemma is a precursor to \lemref{low_density}.

\begin{proof}
  Let $Q_t$ be the number of particles in $S$ that have survived until time
  $t$ but have not left $S$. We have then that $N^n_S(t)\leq Z+Q_t$. The
  rate of coalescence inside $S$ at time $t$ is
  \[
  \sum_{v\in S} \binom{X_v(t)}{2} \geq |S| \binom{Q_t/|S|}{2}
  \]
  (by Jensen's inequality for $\binom{x}{2}$.) If $Q_t<2|S|$ for some
  $t\leq\tau$ then we are done (since $\tau\leq2$.) Otherwise, $|S|\cdot
  \binom{Q_t/|S|}{2} \geq \frac1{2|S|} \binom{Q_t}{2}$, and so $Q_t$ is
  stochastically dominated by the block counting process of a
 Kingman coalescent slowed down by a factor of $2|S|$.
  \lemref{kingman ldev} completes the proof.
\end{proof}

\section{Finite time behaviour of the spatial Kingman coalescent}
\label{S:finitetimeKingman}

\subsection{An induction}

The first step of the argument is to show that for some $m$ (close to
$\log^* n$) there are no particles outside $B(o,m)$ at some specified time,
and to provide lower and upper bounds (both polynomial in the volume of the
ball) on the number of particles at each site inside the ball at the same
time, on an event of high probability. This can be done for any $n$, but it
is easier to consider initially a sub-sequence of $n$'s, and then
interpolate to get the result for all $n$. With this in mind, for a given
integer $m$ denote $V_m = \Vol B(o,m)$, and define $n=n(m):=\Tow(m,V_m^2)$.
Note that $\log^* n = m + \log^* V_m^2$ is very close to $m$, as $m\to
\infty$.

Define the sequence of times $t_k = (\Tow(m-k, V_m^2))^{-3}$, where
$k=0,\ldots,m$. This sequence is increasing from $t_0=n^{-3}$ to
$t_m=V_m^{-6}$. Moreover, $t_{k+1} \gg t_k$ (in particular $t_{k+1} >
2t_k$). Recall that $\theta = 2 \rho$, where $\rho$ is the jump rate of
particles, and that $\Delta$ is the maximal degree in the graph. Define the
events $B_k$ by
\begin{equation}\label{eq:B_k_def}
  B_k = \Big\{ X_v(t_k) = 0 \text{ for all } v\notin B(o,k) \Big\}
  \bigcap
  \left\{ X_v(t_k) \in \left[\tfrac{\rho}{\Delta} t_k^{-1/3}, 4 t_k^{-1}\right]
    \text{ for all } v\in B(o,k) \right\},
\end{equation}
We are particularly interested in the event $B_m$ which states that at time $t_m$
each site of $B(o,m)$ has between $c V_m^2$ and $C V_m^6$ particles, with
no remaining particles outside $V_m$.

\begin{lemma}\label{L:B_m_likely}
  With the above notations, $\P(B_m) \to 1$ as $m\to\infty$.
\end{lemma}

The idea is to prove a bound on $\P(B_k^c)$ by induction on $k$. For $k=0$
we have $\P^n(B_0) \geq 1-2n^{-1}$, since the probability of a pair of
particles coalescing by time $t_0=n^{-3}$ is at most $\binom{n}{2} t_0$,
and the probability of a particle jumping by that time is at most $t_0 n$.
The key to the induction step is the following

\begin{lemma}\label{L:k_step}
  Fix constants $a_0,a_1,\eps>0$. Consider the coalescent started with $n$
  particles, all located at $u\in G$: $X_v(0)=n\delta_{u}(v)$. Let
  $\tau=a(\log n)^{-3}$ for some $a\in[a_0,a_1]$, and define the event
  \[
  A = \cap_v \left\{ X_v(\tau) \in [(1-\eps)Q_v,(1+\eps)Q_v] \right\},
  \]
  where
  \[
  Q_v = \begin{cases}
    2/\tau & v=u, \\ (\theta/d_u) \log n & |v-u|=1, \\ 0 & |v-u|>1.
  \end{cases}
  \]
  Then there exists a $C$ depending on $\eps,a_0,a_1,d_u$ only such that
  \[
  \P^n(A^c) < \frac{C}{\log n}.
  \]
\end{lemma}

\begin{proof}
  In this argument, the expression {\em with high probability (w.h.p.)} stands for
``with probability greater or equal to $1-\frac{C}{\log n}$''.
Let $Z(t)$ be the number of distinct labels corresponding to
  particles that exit $u$ during $[0,t]$ (where each label is counted at most
  once). Let $N(t)$ denote the total number of particles in the coupling
  with the mean-field coalescent of Lemma \ref{L:coupl}. Thus we have:
  \begin{equation}\label{eq:X_NZ2}
    N(t)-Z(t) \leq X_u(t) \leq N(t) + Z(t), \mbox{ almost surely}.
  \end{equation}
Therefore one needs to estimate $N(t)$ and $Z(t)$. For any fixed $\eps$, by
  \lemref{kingman ldev} we have
  \begin{equation}\label{eq:Ntau_conc}
    \P^n(|\tau N(\tau)/2-1| > \eps)
    <  C e^{-c/\tau}
    <  C n^{-1}.
  \end{equation}
So the event $\{\frac{N(\tau)}{2/\tau} \in (1-\eps,1+\eps)\}$ happens with high probability.
 By \lemref{esf ldev} we have
  \begin{equation}\label{eq:Ztau_conc}
    \P^n \left(\left|\frac{Z(\infty)}{\log n} - \theta \right| > \eps
    \right) < C n^{-c}.
  \end{equation}
  In the rest of the argument consider the the process on the event
$B:=\{\frac{N(\tau)}{2/\tau} \in (1-\eps,1+\eps)\} \cap
\{\frac{Z(\infty)}{\log n}\in  (\theta-\eps,\theta+\eps)\}$ that occurs with high probability.
  Note that on $B$, $Z(\tau) \leq Z(\infty)\leq (1+\eps){\log n} \ll 2/\tau$ and
\eqref{eq:X_NZ2} imply the required bounds for $X_u(\tau)$.

Moreover, on $\{Z(\tau)\leq (1+\eps){\log n} \} \supset B$, the probability that at least one particle jumps more than once before time $\tau$ is bounded by $\rho\tau (1+\eps){\log n}=C/\log^2{n}$. On the event that no particle jumps more than once, there cannot be any particle located at
a distance strictly greater than $1$ from $u$ at time $\tau$.

Similarly, on $\{Z(\tau) \leq (1+\eps){\log n} \}\supset B$, the probability of at least one coalescence
  event involving particles located at site $v\neq u$ before time $\tau$ is at most $\tau
  \binom{(1+\eps){\log n}}{2}$, again bounded by $C/\log n$.
 We conclude that w.h.p.~there is no coalescence outside of $u$ before time $\tau$.

  This implies that w.h.p.~the particles located at a neighbor $v$ of $u$ at time $\tau$ are
  precisely those that made a (single) jump from $u$ to $v$.
  To show that their number is close to $(\theta/d_u)\log n$, it suffices to show that
  $Z(\tau)$ is concentrated around $\theta\log n$ (which is already known for
  $Z(\infty)$). Namely, since (on the event of high probability) each jump is to made from $u$
  to a random neighbor of
  $u$, and there are no further moves or coalescence events involving the particles outside of $u$,
  $X_v(\tau)$ is concentrated near $Z(\tau)/d_u$
  for any $v\sim u$, due to a law of large numbers argument.
  Indeed, the number of particles jumping from $u$ to any particular
  of its neighbors has variance of order $\log n$, and using a normal approximation
  to binomial random variables, the probability of deviating by $\eps \log
  n$ from the mean is no more than $C n^{-c}$.

  Thus it remains to show that $Z(\infty)-Z(\tau) \leq \eps\log n$ with high probability.
  To this end, note that $Z(\infty)-Z(\tau)$ is the number of particles
  that exit $u$ after time $\tau$.
  Denote by $\FF_\tau$ the $\sigma$-field generated by the evolution of the process
  up to time $\tau$.
 By \lemref{esf ldev}, monotonicity
  and the Markov property at time $\tau$,
  \begin{align*}
    \P^n(Z(\infty)-Z(\tau) > \eps\log n | \FF_\tau)\indic{N(\tau) < 3/\tau}
    &< \P^{3/\tau} \left(Z(\infty) > \eps \log n \right) \\
    &< C(3/\tau)^C e^{-\eps\log n} \\
    &\ll 1/\log n, \qedhere
  \end{align*}
 and since $\{N(\tau)<3/\tau\}\supset B$
occurs w.h.p., this concludes the argument.\end{proof}

\begin{proof}[Proof of \lemref{B_m_likely}]
  We have $\P(B_m^c) \leq \P(B_m^c \cap B_{m-1}) + \P(B_{m-1}^c)\leq
\P(B_0^c) + \sum_{k\leq m} \P(B_k^c | B_{k-1})$.
  Noting that $\P(B_0^c) \leq 1/n$, we turn to estimating $\P(B_k^c | B_{k-1})$.

Given $\FF_{t_{k-1}}$, consider now the coupling of \lemref{sum_dom}
applied to the process observed on $[t_{k-1},t_k]$,
where the partition $\pi$ is $\FF_{t_{k-1}}$ measurable and where two labels belong
to the same equivalence class of $\pi$ if and only if their corresponding particles
have the same position at time $t_{k-1}$.
On the event $B_{k-1}$ we have that $\log X_u(t_{k-1}) \asymp
  -\log t_{k-1} \asymp t_k^{-1/3}$,\ $u\in B(o,k-1)$, hence \lemref{k_step} applies
to each of the corresponding processes. We
  conclude that with probability at least $1-\frac{C}{\log X_u(t_{k-1})}
  \geq 1 - C t_k^{1/3}$ the following occurs: during $[t_{k-1},t_k]$ ({\em i}) No particle from
  $u$ jumps more than once; ({\em ii}) At most $3/t_k$ particles remain at $u$;
  ({\em iii}) Each neighbor of $u$ receives between $\frac{\rho}{\Delta}
  t_k^{-1/3}$ and $\frac{3\rho}{\Delta} t_k^{-1/3}$ particles. Say that a
  vertex $u\in B(o,k-1)$ is {\em bad at stage $k$} on the complement of the
above event.

  Applying the right hand inequality of \lemref{sum_dom}, on the event $B_{k-1}\cap\{$there are no bad
  vertices at stage $k\}$, we have that $X_u(t_k) = 0$ outside $B(o,k)$ (since
  no particle jumps twice and at time $t_{k-1}$ all the particles are inside $B(o,k-1)$.
  Moreover, each site $v\in B_k$ has at least $\frac{\rho}{\Delta}
  t_k^{-1/3}$ particles jumping to it from some neighbor $u$ of $v$,
  hence $X_v(t_k) \geq \frac{\rho}{\Delta} t_k^{-1/3}$. Finally, for each
  $v\in B(o,k)$ we have $X_v(t_k) \leq 3 t_k^{-1} + 3\rho t_k^{-1/3} < 4
  t_k^{-1}$ (here we may assume that $t_k^{2/3}<1/(3\rho)$), since it receives at
most $\frac{3\rho}{\Delta} t_k^{-1/3}$
  from each of its (at most $\Delta$) neighbors. Hence $B_{k-1}\cap\{$there
is no bad vertex at stage $k\}\subset B_k$.

Therefore $ \P(B_k^c | B_{k-1})\leq \P(\exists$ a bad vertex at stage $k)\leq C
  V_{k-1} t_k^{1/3} \leq \frac{C V_m}{\Tow(m-k, V_m^2)}$.
It follows that
  \[
  \P(B_m^c) \leq \frac{1}{n} + \sum_{k\leq m} \frac{C V_m}{\Tow(m-k,
    V_m^2)} \leq \frac{C}{V_m},
  \]
  since the term for $k=m$ overwhelmingly dominates all the others.
\end{proof}

\subsection{Lower bound estimates}

\lemref{B_m_likely} gives us a fairly accurate description of the spatial
coalescent up to positive times of order $o(1)$.
Additional estimates are needed for understanding the behavior up to a constant time $t$. We
begin with the lower bound, since it is simpler. Henceforth, we let $t>0$
be a fixed time. Recall that the initial configuration of the spatial
coalescent consists of $n=\Tow(m,V_m^2)$ particles located at $o$.

\begin{lemma}\label{L:k-sp1}
  Fix $t>0$. The collection $(X_t(v),\, v\in B(o,m))$
  can be coupled
  with the family $(\zeta_v,\,v\in B(o,m))$ of i.i.d.~Bernoulli variables with mean $e^{-\rho t}$, so that
  \[
  \P\big(\forall v \in B(o,m), X_v(t) \geq \zeta_v \big)
  \xrightarrow[m\to\infty]{} 1.
  \]
\end{lemma}

\begin{proof}
  Assume that $m$ is sufficiently large so that $t_m<t$. By \lemref{B_m_likely}, with
  probability tending to $1$, each site in $B(o,m)$ is not empty at time
  $t_m$. On this event, fix one particle at each $v\in B(o,m)$, and color
  it red. Consider the evolution with coloring (see the proof of \lemref{sum_dom} for a similar construction), so that if a red particle
  coalesces with another particle, the newly formed particle retains the red
  color. Now, it is obvious that between time $t_m$ and $t$, each red
  particle has probability $e^{-\rho (t-t_m)}>e^{-\rho t}$ of not
  migrating, independently of all other red particles, so the claim holds.
\end{proof}

\subsection{Upper bound estimates}

After time $t_m$, the bounds in the definition of $B_m$ (cf.~\eqref{eq:B_k_def})
still hold for most vertices, but will begin to fail for some vertices. As
the number of particles per vertex decreases, the probability of failure
increases. We overcome this by combining the second part of \lemref{esf
  ldev} with \lemref{k-few in S}.


\begin{lemma}\label{L:k-sp2}
  Fix $\eps,t>0$, and start with $n=\Tow(m,V_m^2)$ particles at $o$. With
  high probability there is no particle outside $B(o,(1+\eps)m)$ at or
  before time $t$, and the total number of particles at time $t$ is at most
  $C V_{(1+\eps) m}$.
\end{lemma}

\begin{proof}
  By \lemref{B_m_likely}, with high probability at time $t_m$ there are no
  particles outside $B_m$, and the number of particles inside $B(o,m)$ is
  at most $4 V_m/t_m = 4 V_m^7$. By ignoring coalescence transition after
  time $t_m$, so that each particle performs a simple random walk
  independently of all the others, the number of particles located at any
  particular site at any later time can only become larger. Each particle
  makes an additional $\Poi(\rho (t-t_m))$ steps during $[t_m,t]$, so the
  probability that at least one of these particles makes at least $\eps m$
  steps is bounded by $4 V_m^7 C^{\eps m} / \lfloor \eps m \rfloor!$. This
  last quantity tends to $0$, since $V_m \leq C \Delta^m$.

  Thus with an overwhelming probability, there are no particles outside
  $B(o,(1+\eps)m)$ at time $t$. By \lemref{k-few in S} and the above
  observation, the number of particles within $B(o,(1+\eps)m)$ is at most a
  constant multiple of $V_{(1+\eps)m}$, again with an overwhelming
  probability.
\end{proof}

\subsection{Interpolation}

\begin{proof}[Proof of \thmref{k finite}]
  If $n=n(m)=\Tow(m,V_m^2)$ for some $m$, then Lemmas~\ref{L:k-sp1} and
  \ref{L:k-sp2} imply that with high probability the number of particles at
  time $t$ is between $c V_m$ and $C V_{(1+\eps)m}$. Since $\log^* n = m +
  \log^* V_m^2 \sim m$, this implies the claim.

  For intermediate $n$, we use the monotonicity of the process in $n$. Note
  that $\log^* n(m+1) - \log^* n(m) \leq 2$ (since $V_{m+1} < e^{V_m}$), so
  that the sequence $n(m)$ is sufficiently dense to imply the theorem.
\end{proof}

\begin{rem*}
  Since $m = \log^* n(m) - \log^* m + O(1)$, the proof above gives the
  lower bound $c V_{m-\log^* m}$. As for the upper bound, the proof of
  \lemref{k-sp2} works with radius $m + C m/\log m$ in general, and $m +
  \log m$ for graphs with polynomial growth.

  It is possible to get both lower and upper bounds that are closer to
  $\Vol B(o,\log^* n)$. For the lower bound, one way would be to argue that
  most vertices continue to behave typically (as in \lemref{k_step}) even
  up to constant times.

  The upper bound is more delicate. One way of improving it is by
  considering the evolution of the total number of particles in $B(o,k)$
  for $k>m$, similarly to the argument of Section~\ref{sec:long time}. Under additional growth
  assumptions on the graph, both bounds are of order $\Vol B(o,\log^*
  n)$.
\end{rem*}

\section{Results for spatial Beta-coalescents}
\label{S:ProofBeta}

We now turn to the proof of \thmref{beta}. In fact, we prove a slightly
more general result. Suppose that $\Lambda$ has a  sufficiently regular density near $0$:
$\Lambda(d x)=g(x) dx$, where for some $B>0$ and
$\alpha<\in (1,2)$ we have
\begin{equation}\label{reg-var}
  g(x)\sim B x^{1-\alpha}, \qquad x\to 0.
\end{equation}
This includes the case where $\Lambda$ is the Beta$(2-\alpha,\alpha)$
distribution. A consequence of \eqref{reg-var} is the following standard
estimate for the rate of coalescence events when there are $n$ particles
remaining:

\begin{lemma}\label{L:lambda_asymp}
 The sequence $(\lambda_n)_{n\geq 2}$ is increasing in $n$. Furthermore, there exists $c>0$ which
  depends only on $\alpha,B$, such that if $\Lambda$ satisfies
  \eqref{reg-var}, then $\lambda_n \sim c n^{\alpha}$.
\end{lemma}

\begin{proof}
  The monotonicity of $\lambda_n$ in $n$ is a consequence of the natural
  consistency of $\Lambda$-coalescents. The second part of the statement is
  a consequence of \eqref{reg-var} and Tauberian theorems. See, e.g.,
  \cite[Lemma~4]{blg3} for more details.
\end{proof}

\subsection{Lower bound in Theorem \thmref{beta}}

Define the following parameters
\begin{align}\label{gamma}
  \beta &= \frac{\alpha-1}{2} & \tau &= an^{-\beta} \text{ for some } a \in
  [a_0,a_1]& \gamma &= \min\{1-\alpha/2,\beta/2,1/8\},
\end{align}
and observe that both $\gamma>0$ and $\alpha-2+\gamma \le -\gamma$. We
next consider the quantity
\[
Y_n= \int_0^\tau N(s) \ ds.
\]

\begin{lemma}\label{Yn}
  Assume that $\Lambda$ satisfies \eqref{reg-var}. Then for some $c,C$
  depending only on $\Lambda$,
  \begin{equation}\label{Xn-ub}
    \P(Y_n \ge n^{2-\alpha+\gamma})\le C n^{-\gamma}, \ \forall n\geq 2,
  \end{equation}
  and
  \begin{equation}\label{Xn-lb}
    \P(Y_n \le cn^{\gamma}) \le Cn^{-\gamma}, \ \forall n\geq 2.
  \end{equation}
\end{lemma}

\begin{rem*}
  It follows from Theorem 5 in \cite{bbl1} that $Y_n \sim c n^{2-\alpha}$,
  almost surely as $n\to \infty$, for some $c>0$. However this result does
  not provide any estimate on the deviation probability.
\end{rem*}

\begin{proof}
  The key fact is that if the process $N(t)$ attains some value $k$, then
  it stays at $k$ for an exponentially distributed time with mean
  $1/\lambda_k$. Since the probability of hitting $k$ is at most 1,
  \[
  \E Y_n \leq \sum_{k\leq n} \frac{k}{\lambda_k} \leq cn^{2-\alpha}
  \]
  by \lemref{lambda_asymp}. The upper bound \eqref{Xn-ub} follows by
  Markov's inequality.

  The lower bound is more delicate. We argue that with high probability the
  first $M=n^{\alpha-1+\gamma}$ jumps all occur before time $\tau$ and that
  throughout these jumps $N(t)$ remains above $n/2$. Summing over only these
 jumps will give the lower bound \eqref{Xn-lb}.

  Let $B_m$ be the number of particles lost in the next coalescence when
  there are $m$ particles present. It is known \cite[Lemma~7.1]{bbs1} that
  there exists $C>0$ such that
  \begin{equation}\label{jump sizes}
    \P(B_m>k) \le Ck^{-\alpha} \text{ for all }m,k\ge 1.
  \end{equation}
  In particular, $\E B_m < c$ for some constant depending only on
  $\Lambda$. Thus the total size of the first $M$ jumps has expectation at
  most $cM$. Let $t_k$ be the time of the $k$th jump in $N(t)$, then by
  Markov's inequality
  \begin{equation}\label{E:ubjumps}
    \P^n(N(t_M)<n/2) < \frac{cM}{n-n/2} < cn^{\alpha-2+\gamma} < cn^{-\gamma}.
  \end{equation}

  On the event that $N(t_M) \geq n/2$, the rate of each of the first $M$
  jumps is at least $\lambda_{n/2}$. Thus, by Markov's inequality, and by
  monotonicity of $\lambda_m$,
  \begin{equation}\label{E:ubjumps2}
    \P(t_M > \tau, N(t_M)\geq n/2) \leq \frac{M/\lambda_{n/2}}{\tau}
    \leq cn^{-1+\gamma+\beta} < cn^{-\gamma}.
  \end{equation}
  Thus, combining (\ref{E:ubjumps2}) with (\ref{E:ubjumps}),
  $\P(A^c)<cn^{-\gamma}$, where $A=\{t_M < \tau, N(t_M)\geq n/2\}$.

  Note that, on the event $A$,
  \[
  Y_n = \int_0^\tau N(t) dt \geq \int_0^{t_M} N(t) dt \ge (n/2)t_M.
  \]
  It thus suffices to show that $\P^n(t_M \le cn^{\gamma-1}) \le
  Cn^{-\gamma}$. However, the rate of each jump is at most $\lambda_n$, and
  therefore
  \[
  t_M \sgt \sum_{i=1 }^M E_i
  \]
  where $E_i$ are i.i.d.\ exponentials with rate $\lambda_n$. Now, from
  \lemref{lambda_asymp} we know that
  \[
  \E \sum_{i\leq M} E_i \sim c n^{\gamma-1},
  \]
  and by Lemma \ref{L:ldev} with $\epsilon =1/2$,
  \[
  \P\left(\sum_{i\leq M} E_i < c n^{\gamma-1}/2\right) <
  \exp\left(-\frac1{16}n^{\alpha-1+\gamma}\right) < Cn^{-\gamma}
  \]
  as needed. This completes the proof of Lemma \ref{Yn}.
\end{proof}

The next result gives a lower bound on the number of particles that exit
the origin. This complements the upper bound of \lemref{few_leave}. Recall
that $Z(t)$ is the number of particles that exit the origin by time $t$.
The idea is that as long as $Z$ is small, the true behavior is close to the
upper bound.

\begin{lemma}\label{L:many_leave}
  Let $A$ be the event $\{Z(\tau) < n^\gamma\}$. Then $\P(A) =
  O(n^{-\gamma})$.
\end{lemma}

\begin{proof}
  We introduce the random time $T_a$ defined for any $0<a<1$ by $T_a=
  \inf\{t>0: Z(t) \ge aN(t)\}$. Define
  \begin{align*}
    A_1 &= \{Z(\tau\wedge T_a) \le n^\gamma\} & A_2&= A \cap \{ \tau>T_a\}.
  \end{align*}
  Note that $A\subset A_1\cup A_2$ so it suffices to prove that
  $P(A_i)=O(n^{-\gamma})$, for $i=1,2$.

  Consider $A_1$ first. Recall the notations introduced in Lemma
  \ref{L:few_leave}, and note that $T_a$ is a stopping time with respect to
  the filtration $\FF^*$. Since $N(t)$ is non-increasing with limit 1 and
  since $Z(t)$ is non-decreasing and non-negative integer valued, $T_a$ is
  finite if and only if at least one particle leaves $o$. This will
  eventually happen, so $T_a$ is a.s.\ finite. Denote by $\tilde \P_n$ the
  law
  \[
  \tilde \P_n(\cdot)=\P_n(\cdot|\GG),
  \]
  of all processes, conditioned on the entire evolution of $N$.

  Consider the martingale $S_t$ stopped at time $T_a$. By Doob's
  inequality, we find that for any $\delta>0$
  \begin{align}
    \tilde \P_n \left( \sup_{s\leq T_a} |S_s| \geq \delta \int_0^{T_a} N(s)
      ds \right) &\leq \frac{4\rho \tilde\E_n\left(\int_0^{T_a} M_u
        du\right)}
    {\delta^2\left(\int_0^{T_a} N(s) ds\right)^2}\wedge 1\nonumber\\
    &\leq \frac{4\rho}{\delta^2\int_0^{T_a} N(s) ds} \wedge 1.
    \label{doob}
  \end{align}
  The last inequality follows from the first bound of \eqref{EMZN}, which
  implies that $\tilde \E_n (\int_0^{T_a} M(u) du )\le \int_0^{T_a} N(u)
  du$. Define the event
  \[
  A_s = \left\{1-a-\delta < \frac{Z_s}{\rho\int_0^s N(u) \,du} < 1+\delta
  \right\}.
  \]
  Until time $T_a$ we have $M(t)\ge(1-a)N(t)$, and so \eqref{EMZN} and
  \eqref{doob} imply
  \[
  \tilde\P_n(A_s^c) \leq \frac{4\rho}{\delta^2\int_0^{T_a} N(s) \,ds}.
  \]
  We fix $a$ and $\delta$ such that $1-a-\delta>1/2$. After taking the
  expectation, we obtain, using \eqref{Xn-ub}:
  \[
  \P^n(A_1) \le O(n^{-\gamma}) + n^{\alpha-2-\gamma} = O(n^{-\gamma})
  \]

  Turning to $A_2$, note that
  \[
  A_2 \subset \{a N(\tau) \le n^\gamma\}
  \]
  We claim that
  \begin{equation}\label{Ntau}
    \P^n(a N(\tau) \le n^\gamma) \le Cn^{-\gamma}.
  \end{equation}
  To see this, we use the following rough estimate. Note that by
  \eqref{jump sizes}, there is a probability at least
  $1-Cn^{-\gamma\alpha}$ that $N(s) \in [n^\gamma+1,2n^\gamma]$ for some
  $s$. In this case, the process will wait an amount of time greater than
  an exponential $Y$ with rate $\lambda_{2n^\gamma}$ before the next jump.
  It follows that (since $\gamma\le \beta/2$ and $\alpha<2$),
  \begin{align*}
    \P^n(N(\tau) \le n^\gamma)
    &\le \big( 1 - C n^{-\gamma \alpha} \big) \P(Y\le \tau) \\
    &\leq 1 - C n^{-\gamma \alpha} - \exp(-c\tau n^{\alpha\gamma}) \\
    &\leq 1 - C n^{-\gamma \alpha} - \exp(-cn^{-\gamma}) \\
    &< cn^{-\gamma}.
  \end{align*}
  This completes the proof of \lemref{many_leave}.
\end{proof}

We are now ready to start proving the lower-bound of \thmref{beta}. Let
$t_0>0$ be a fixed time.

\begin{lemma}\label{L:B_step-lb}
  Fix constants $a_0,a_1$ such that $1<a_0<a_1$. Consider the coalescent started with $n$
  particles, all located at $u\in G$: $X_v(0)=n\delta_{u}(v)$. Let
  $\tau=an^{-\beta}$ for some $a\in[a_0,a_1]$, and define the event $A$ by
  \[
  A=\{X_u(\tau) \ge n^\gamma/(4d_u)\} \cap \cap_{v\sim u}\{X_v(\tau) \ge n^\gamma/(4d_u)\}.
  \]
  There are constants $c,C$ depending on $a_0,a_1,d_u$ only such that
  $\P(A^c) < C n^{-c}$.
\end{lemma}

\begin{proof}
  The fact $\P(X_u(\tau) < n^\gamma/(4d_u))<C n^{-c}$ is a direct
  consequence of  \eqref{Ntau} where we choose $a<1<4d_u$ satisfying
  $1-a-\delta>1/2$.
  For $v\sim u$, \lemref{many_leave} gives a bound on the probability that
  not many particles leave the origin. It is highly probable that
  a proportion close to $1/d_u$ of these particles jumps to $v$. It remains to
  estimate the number of particles that move to $v$ and subsequently
  coalesce.

  If all the particles that migrate to $v$ do so immediately at time 0,
  so that they have strictly more opportunities to coalesce, the number of particles
  remaining at $v$ at time $\tau$ would still be sufficiently large.
  Indeed, it would then take $Y$ amount of time, where $Y$ is an exponential random
  variable with parameter
  $\lambda_{n^\gamma/(4d)}$, before the first coalescence. Since
  $\lambda_m \le cm^\alpha$ for all $m\ge 1$, we deduce that $\E(Y)\ge
  cn^{-\gamma\alpha}$. However, since $\gamma\le \beta/2$ and $\alpha<2$, we have
   $\tau=an^{-\beta}\ll cn^{-\gamma\alpha}$, hence
  $\P(Y<\tau)\le cn^{\alpha\gamma - \beta}$.

  In addition, note that by Lemma \ref{L:lambda_asymp}, the total jump rate of
  $n^\gamma$ particles is smaller than the total coalescence rate (since
  $\alpha >1$), so the probability any of the particles that jump to $v$ makes an
  extra jump before time $\tau$ is smaller than $cn^{\alpha\gamma - \beta}$.
 It follows that  there are at least $n^\gamma/(4d_u)$ particles located at $v$
 at time $\tau$, with
  probability greater than $1-Cn^{-c}$.
\end{proof}

\begin{proof}[Proof of \thmref{beta} (lower bound)]
  Let $f_k(n)=f\circ f\ldots \circ f(n)$ ($k$ iterations) where $f(n)=
  n^\gamma/4d$. Define the sequence of times $(\tau_k)_{k=1}^\infty$
  \[
  \tau_k= \tau_{k-1}+af_{k-1}(n)^{-\beta}.
  \]

  It is easy to check that if we take $k=k(n)= \log \log n/(-2\log
  \gamma)$, then
  \[
  f_k(n) \ge c \exp(\sqrt{\log n})
  \]
  Let $A'$ be the event that at each site within radius $k$ there are at
  least $f_k(n)$ particles at time $\tau_k$. On $A'$, reasoning as in Lemma
  \ref{L:k-sp2}, (at each site of this ball at least one particle
  remains with positive probability until time $t_0$), we see that
  $N^n(\tau) \ge \Vol B(o,k) \ge c \Vol B(o, c\log \log n)$ for some $c>0$.
  Thus to obtain the lower bound of \thmref{beta}, it suffices to compute
  the cumulative error probability in the iterated application of
  \lemref{B_step-lb}. However, it is easy to check that
  \[
  \P(A'^c) \le \sum_{i=1}^{k} C\Vol B(o,i) f_i(n)^{-\gamma} \le C k \Vol
  B(o,k) f_k(n)^{-\gamma}.
  \]
  Since $\Vol B(o,k) < \Delta^k$, where $\Delta$ is the degree of the graph,
  this converges to 0 as $n\to\infty$.
\end{proof}

\subsection{Upper bound in Theorem \thmref{beta}}

The proof of the upper-bound in Theorem \thmref{beta} requires a few
additional estimates.

Consider the spatial coalescent on any graph $G$.
  Given some  subset $A\subset V$ of the vertices,
  denote by $Q_t$ the number of particles that are present in $A$
  throughout the time interval $[0,t]$.
\begin{lemma}\label{L:B_few_remain}
There are constants $c,C>0$ which depend  on $\Lambda$ only, so that
  \[
  \P(Q_{t_0}>C t_0^{-1/(\alpha-1)} |A|)<\exp(-c|A|).
  \]
\end{lemma}

\begin{proof}
  Ignoring the particles after they exit $A$,
  one may assume that any particle leaving $A$ is immediately killed.
  The main reason for $Q_{t_0}$ being small is the coalescence.
  The total rate of coalescence at a site $v$ holding $X_v$ particles is
  $\lambda_{X_v}\sim c (X_v)^\alpha$. At each such event at least one
  particle disappears, and therefore the total rate of decrease of $Q_t$ at time $t\le t_0$ is
  at least
  \[
  \sum_{v\in A} (c X_v(t))^\alpha \geq c |A|^{1-\alpha} Q_t^\alpha,
  \]
  due to Jensen's inequality, since $\alpha >1$. (This
  is similar to \cite[Theorem~12]{LS}, but the above inequality is
  stronger). Thus $(Q_t,t\le t_0)$ is stochastically dominated by a pure death chain where the rate
  of decrease from $i$ to $i-1$ is $c |A|^{1-\alpha} i^\alpha$.

 One concludes the argument using \lemref{ldev}. Let $E_k$ be independent
  exponential random variables with mean $\mu_k=c |A|^{\alpha-1}
  k^{-\alpha}$, and define $S_K = \sum_{k>K} E_k$. Then we have
  \[
  \P(Q_{t_0}>K) < \P(S_K>t_0).
  \]
  To apply \lemref{ldev} to $S_K$ we need to estimate $\E S_K $ and $ \Var
  S_K$: note that for suitable constants, as $K\to\infty$,
  \begin{equation} \label{exp-ub} \E S_K = \sum_{k>K} \mu_k^{-1} \sim c_1
    |A|^{\alpha-1}K^{1-\alpha}
  \end{equation}
  and
  \begin{equation}
    \Var S_K =\sum_{k>K} \mu_k^{-2} \sim c_2 |A|^{2\alpha-2} K^{1-2\alpha}.
  \end{equation}
  In particular $\frac{\Var S_K}{\E S_k \mu_K}$ is asymptotically constant
  and we may apply \lemref{ldev} with some constant $\eps$. Thus for some
  $c_3>0$,
  \begin{align*}
    \P(S_K > 2\E S_K) &\le \exp\left(- c \frac{ (\E S_K)^2}{\Var S_K}
    \right) \\
    &< e^{- c_3 K }.
  \end{align*}

  Now, if $K$ is such that $\E S_K < t_0/2$ we may conclude that
  \[
  \P(Q_{t_0}>K) < e^{-c_3 K}.
  \]
  From \eqref{exp-ub} we see that $K=Ct_0^{-1/(\alpha-1)}|A|$ works for $C$
  large enough.
\end{proof}

\begin{lemma}\label{L:B_step-ub}
  Fix constants $a_0,a_1,\eps>0$. Consider the coalescent started with $n$
  particles, all located at $u\in G$: $X_v(0)=n\delta_{u}(v)$. Let
  $\tau=an^{-\beta}$ for some $a\in[a_0,a_1]$, and define the event $A$ by
  \[
  A= \bigcap_v \{X_v(\tau) \le C_1 Q_v\},
  \]
  with
  \[
  Q_v = \begin{cases}
    n^{3/4} & \text{ if } v=u, \\
    n^{2-\alpha+\gamma} & \text{ if } |v-u|\le  r:=\lceil 4/(\alpha-1)\rceil,\\
    0 & \text{otherwise}.
  \end{cases}
  \]
  Then there are constants $C,C_1$ depending only on $\Lambda,a_0,a_1$ such
  that $P(A^c) < C n^{-\gamma}$.
\end{lemma}

\begin{proof}
  With sufficiently high probability at most $n^{2-\alpha+\gamma}$
  particles leave the origin by time $\tau$ (due to \lemref{few_leave} and
  \eqref{Xn-ub}). This implies the bound for $0<|v-u| \le r$.

  Some of the at most $n^{\alpha -2 +\gamma}$
particles leaving $u$ may coalesce before time
  $\tau$, but this may only reduce further the number of particles. We claim that
  except on an event of polynomially small probability, none of these particles makes more
  than $r$ jumps by time $\tau$. Indeed, the probability that by time
  $\tau$, a given particle has jumped more than $r$ times is smaller than
  $(\rho a n^{-\beta})^r$ and there can never be more than $n$
particles in total. Thus if $r$ is such that $n^{1-r\beta}<n^{-\gamma}$,
the  probability of any particle reaching distance $r$ is indeed smaller than
  $Cn^{-\gamma}$, implying the statement of the lemma for any $v$ such that $|v-u|\geq r$.
   For the case $v=u$ we
  invoke \lemref{B_few_remain} with an arbitrary set $A\ni u$ of
  size $c\log n$. If $c$ is large enough then, except on an event of
  probability bounded by
  $n^{-\gamma}$, we have $Q_\tau < C n^{1/2}|A| \ll n^{3/4}$. However,
  $X_u(\tau) < Q_\tau+Z_\tau$, and so by \lemref{few_leave} and
  \eqref{Xn-ub} again, $X_u(\tau) \ll n^{3/4}$. It is easy to see from
  \eqref{gamma} that for all $\alpha\in(1,2)$ we have $2-\alpha+\gamma <
  3/4$.)
\end{proof}

\begin{proof}[Proof of \thmref{beta}: upper bound]
  Note that for any $\alpha$ we have $\gamma\leq\beta/2<1/2$. Let
  $c=\max(\gamma,3/4)$, and note that $c<1$. Let $C_2=C_1\times \Vol B(o,
  r)$, where $C_1$ and $r$ are the constants in \lemref{B_step-ub}.

  Let $f(n)=C_2 n^{3/4}$, and as before set $f_k(n) = f\circ \cdots \circ
  f(n)$ ($k$ iterations). Also set $\tau_1=\tau=n^{-\beta}$, and
  \[
  \tau_{k}= \tau_{k-1}+ a f_{k-1}(n)^{-\beta}.
  \]
  Let $A_i$ be the event that at time $\tau_i$ there are no particles
  outside $B(o,ir)$ intersected with
  \[
  \bigcap_{v,\, |v|\le ir}\{X_v(\tau_i) \le f_i(n)\}.
  \]

  Choose $k=k(n)$ to be the maximal $k$ so that $f_k(n)>\log n$. It is
  clear that $f_k(n) < (\log n)^2$. It is also straightforward to check
  that $k\sim c\log\log n$, and that $\tau_k =o(1)$.

  Applying \lemref{B_step-ub} iteratively, we see that
  \begin{equation}\label{error-ub}
    \P(A_k^c) \le \sum_{i<k} Cf_i(n)^{-\gamma}\Vol(B(o,ir))
    \le C \Vol(B(o,kr)) f_k(n)^{-\gamma} \xrightarrow[n\to \infty]{} 0.
  \end{equation}
  Consequently, at time $\tau_k$ the total number of
  remaining particles is at most $Cf_k(n)\Vol(B(o,kr))$, and these particles are
  all located in $B(o,kr)$, with high probability.

  Consider now the set $B'=B(o,M\log\log n)$ for some large $M$ to be
  specified soon. In order for any particle to exit $B'$ by time $t$ it must
  survive to time $\tau_k$ and jump at least $M\log\log n - kr$ times by
  time $t$. Thus the expected number of particles that exit $B'$ by time $t$
  is at most
  \[
  Cf_k(n)\Vol(B(o,kr)) e^{-c(M\log\log n - kr)} < C(\log n)^2 (\log\log
  n)^d e^{-(cM-c')\log\log n}.
  \]
  Fix $M$ large enough that the last expression tends to 0 as $n\to \infty$.

  Finally note that if no particle leaves $B'$ then $\sum_v X_v(t) = Q_t$. By
  \lemref{B_few_remain}, with high probability the number of particles that
  remain in $B'$ throughout $[0,t]$ is at most $O((\log\log n)^d)$.
\end{proof}

\section{Global divergence of spatial $\Lambda$-coalescents}
\label{S:global}

\subsection{Infinite tree length for $\Lambda$-coalescents}

Fix an arbitrary probability measure $\Lambda$ on $[0,1]$. Consider the
corresponding mean-field $\Lambda$-coalescent that starts from a
configuration consisting of infinitely many blocks, and let $(K^n(s),s\ge
0)$ be the number of blocks process of its restriction to the first $n$
particles.
Define:
\begin{equation}\label{Dint}
  X_n(t)\equiv X_n = \int_0^t (K^n(s)-1) ds.
\end{equation}
The notation $K^n$ might be suggestive of the Kingman coalescent,
so we wish to point out that the measure $\La$ in the following calculation is quite general.

We are interested in the quantity $X_n$ due to the following observation:
if $K^n$ is a good approximation for the number of blocks at the origin of
the spatial $\Lambda$-coalescent at small times $s$, then for $t$ small,
$\rho X_n$ approximates well the number of particles that emigrate from the
origin up to time $t$ (see, for instance, Lemma \ref{L:few_leave}). The key
ingredient in the proof of \thmref{global} is the following result.

\begin{lemma}\label{L:intdiv}
  For any fixed $t>0$ we have $X_n\xrightarrow[n\to\infty]{} \infty$ almost
  surely.
\end{lemma}

\begin{proof}
  Denote by $\sim^t$ the equivalence relation on the labels generated by
  the coalescent blocks at time $t$. For $n\geq 2$ let
  \[
  \tau_n := \min\{t > 0: \exists j < n \mbox{ s.t. } n\sim^t j\}
  \]
  be the first time that the particle labelled $n$ coalesces with any of
  the particles with smaller labels. We have that
  \[
  K^n(s) = K^{n-1}(s) + \indic{s<\tau_n},
  \]
  and therefore
  \[
  X_n = X_{n-1} + (\tau_n\wedge t),
  \]
  i.e.\ the contribution to $X_n$ of particle $n$ is $\tau_n\wedge t$.

  Define $\FF_n$ to be the $\sigma$-algebra generated by $\{K^j_s\}_{j\leq
    n, s>0}$. Conditioned on $\FF_{n-1}$, the infinitesimal rate of
  coalescence of particle $n$ with particles with smaller labels at time
  $s$ is given by
  \[
  \int_{[0,1]} \frac{1}{x^2}\cdot x \cdot (1-(1-x)^{K^{n-1}(s)}) \,
  d\Lambda(x).
  \]
  Applying $(1-x)^k \ge 1-kx$ (for $x\in[0,1]$) we find that the rate of
  coalescence of particle $n$ is at most $K^{n-1}(s)$ (with equality if and
  only if $\Lambda$ is the point mass at 0, in which case the coalescent is
  Kingman's coalescent). Thus
  \begin{align*}
    \E(\tau_n\wedge t | \FF_{n-1})
    &= \int_0^t \P(\tau_n > s|\FF_{n-1}) \, ds \\
    &\ge \int_0^t \exp \left(-\int_0^s K^{n-1}(u) \,du \right)\, ds\\
    &\ge \int_0^t \exp \left(-s-\int_0^t (K^{n-1}(u)-1) \,du \right) \, ds\\
    &= e^{-X_{n-1}} \int_0^t e^{-s} ds = (1-e^{-t}) e^{-X_{n-1}}.
  \end{align*}

  Note that $X_n$ is increasing and consider the martingale
  \[
  M_n = X_n - \sum_{k=2}^n \E(\tau_k\wedge t|\FF_{k-1}).
  \]
  On the event that $X_n$ is bounded, the last calculation implies that
  $\E(\tau_k\wedge t | \FF_{k-1})$ is bounded from below, hence
  $M_n\to-\infty$. Since $M$ is a martingale, the last event has
  probability $0$.
\end{proof}

Note that a different proof of Lemma \ref{L:intdiv} follows from Corollary
3 in \cite{bbl1}, although the arguments there are significantly more
involved.

\subsection{Proof of Theorem \ref{T:global}}

We now consider the spatial coalescent corresponding to some fixed
$\Lambda$ as in the previous section, on an arbitrary locally finite graph
$G$. As usual, let $n$ denote the initial size of the population, with all
particles initially located at $o$, a fixed vertex of $G$. Recall the
definitions of the processes $M$ and $Z$ in Lemma \ref{L:coupl}. Both
processes $M$ and $Z$ depend implicitly on $n$, omitted from the notation.
We consider the usual coupling of coalescents that correspond to different
$n$.

\begin{lemma}
  For any $t>0$ we have that $Z(t)\xrightarrow{}\infty$ almost surely as
  $n\to \infty$.
\end{lemma}

\begin{proof}
  We follow the argument of \lemref{many_leave}, except that we are only
  interested in showing that $Z$ diverges, which simplifies the argument.
  Since $Z(t)$ is non-decreasing in $n$ it suffices to show that for any
  fixed $m$ we have $\P(Z(t)<m)\xrightarrow[n\to\infty]{} 0$.

  Recall the martingales \eqref{E mart S}. On the event $\{Z(t)\le m\}$, we
  have for all $s\le t$ that $M_s\ge N(s)-Z_s \ge N(s)-m$, and therefore
  $S_t \le m + \rho mt - \int_0^t \rho N(t) dt$. Due to \lemref{intdiv},
  for any fixed $m,t$ and any sufficiently large $n$, on the event
  $\{Z(t)\le m\}$ (this event also depends on $n$)
  \[
  S_t \le -\frac12 \int_0^t \rho N(t).
  \]
  As in \eqref{doob}, Doob's maximal inequality yields that for large
  enough $n$
   \[
   \P(Z_t<m|\GG) \le \P \left(\left. \sup_{s\leq t} |S_s| \ge \frac12
       \int_0^{t} \rho N_s ds \,\right|\,\GG\right) \le \frac{16}{\rho
     \int_0^t N_s ds}\,.
  \]
  By \lemref{intdiv} the right-hand side tends to 0 almost surely, so
  $\P(Z(t)<m)\xrightarrow[n\to\infty]{} 0$.
\end{proof}

Fix $\eps >0$ and a vertex $v$ of the graph, and let
$E_{m,\eps,v}=\{\sup_{t\in [0,\eps]} X_v^n(t) \geq m\}$ be the event that
at some time $t<\eps$ there are at least $m$ particles located at site $v$.

\begin{lemma}\label{L:v_reached}
  We have $\P^n(E_{m,\eps,v})\xrightarrow[n\to\infty]{} 1$.
\end{lemma}

\begin{proof}
  Note that the claim is trivially true if $v=o$. We prove it first for $v$
  a neighbor of $o$. Take
  \[
  t_0=\eta\min\{\eps,\lambda_m^{-1},(\rho m)^{-1} \},
  \]
  where $\eta$ is an arbitrarily small number. Now, choose $n_0=n_0(t_0)$
  large enough that $Z_{t_0}>2dm$ with probability at least $1-\eta$, where
  $d=\deg(o)$. By the weak law of large numbers, one can choose $n_0 =
  n_0(t_0,h)$ large enough that on the event $\{Z_{t_0} >2dm\}$, $v$
  receives at least $m$ particles from $o$ with probability at least
  $1-\eta$. We concentrate on this event of high probability, and on these
  $m$ particles, ignoring any further particles that might visit $v$.

  Jumps from $v$ occur at rate $\rho$ per particle, so the probability that
  any of the $m$ particles above leave $v$ before time $t_0$ is at most
  $\eta$. Since a coalescent event involving any $k$-tuple of particles
  occurs at total rate $\lambda_k$ (increasing in $k$), and since at a
  given time $s<t_0$ there are up to $m$ of the above particles located at
  $v$, the probability of a coalescent event before time $t_0$ in which two
  or more of the $m$ above particles participate is at most $\eta$. It
  follows that $\P(X^n_v(t) \ge m) \ge 1-4\eta$, for all $n>n_0$. Since
  $\eta$ can be made arbitrarily small, this proves our claim for $v$ a
  neighbor of $o$.

  For other $v$ we use induction in the distance $|v|$ to $o$. Indeed, such
  $v$ has a neighbor $u$ satisfying $|u| < |v|$. For any fixed $m', \eta$,
  and $n$ sufficiently large, we have $\P^n(E_{m',\eps,u})\geq 1-\eta$.
  Given this, and using the strong Markov property, one can repeat the
  previous argument with $m'$ sufficiently large to conclude that with
  probability at least $1-2\eta$ there will be at least $m$ particles at
  $v$ (arriving from $u$) at some time $t<2\eps$.
\end{proof}

\begin{proof}[Proof of \thmref{global}]
  Again, due to monotonicity in $n$ and $t$, it suffices to show that for
  any $m <\infty$ and any $t>0$, we have $\lim_{n\to \infty} \P(N^n_t >
  m)=1$.

  Let $\eta,\eps >0$ be small numbers. Fix $m<\infty$, and choose a
  subgraph $G^m \subset G$ of size $m$ such that the distance between any
  two vertices of $G^m$ is larger than $1/\eta$. By \lemref{v_reached} we
  have that
  \[
  \P^n(E_{1,\eps,v})\xrightarrow[n\to\infty]{} 1, \ \forall v \in G^m.
  \]
  Moreover, if $A_{v,\eps}$ is the event that the first (if any) particle
  that enters $v$ before time $\eps$ stays at $v$ up to time $\eps$ (while
  it may possibly coalesce with other particles), note that $\P(A_{v,\eps}|
  E_{1,\eps,v}) \geq e^{-\rho\eps}$. By choosing $\eps$ sufficiently small
  we arrive to
  \[
  \lim_{n\to \infty}\P(\cap_{v\in G^m} A_{v,\eps}) \geq 1-\eta.
  \]
  However, given $\cap_{v\in G^m} A_{v,\eps}$, the probability that any
  pair of the above particles (located at mutual distance greater than
  $1/\eta$ at time $\eps$) will coalesce before time $t$ tends to 0 as
  $\eta \to 0$.
\end{proof}

\section{Lower bound for the long time asymptotics}
\label{sec:long time}

We now turn to the large time asymptotic behavior of spatial coalescents.
The underlying measure $\Lambda$ does not play an important role here as it
did for the behavior at constant times. The reason for this is that, as
explained in the introduction (Section \ref{S:ideas}), at the beginning of
this phase, say at constant time $t>0$, the number of particles at each
site is tight with respect to $n$. When the number of particles at a site
is small, the coalescents corresponding to different choices of $\Lambda$
behave similarly. In fact, the density of particles quickly decays, and
once it is small enough, it rarely happens that more than two particles are
at the same location. With at most two particles at each site, any spatial
$\Lambda$-coalescent is equivalent to spatial Kingman's coalescent.

An important quantity in this setting is the radius $m$ of the region
(ball) which is initially ``filled'' with particles. As we have seen, for
Kingman's coalescent the radius of this ball is $m=\log^* n$, while in the
case of Beta-coalescents it is approximately $m =\log \log n$ up to
constants. In the general case, the radius $m$ should be a certain function
of both $n$ and $\Lambda$, namely $m= f^*(n)$ where $f^*$ is defined in
(\ref{fstar}). This was rigorously established only for Kingman's
coalescent and those with ``regular variation'' (i.e., satisfying
(\ref{reg-var})). However, the results which we present in this section and
the next one, are valid for essentially arbitrary coalescence mechanisms
(subject to (\ref{condUB}) for the upper bound in Section \ref{SS:upper
  bound}), and assume that the spatial $\Lambda$-coalescent starts with a
possibly random but tight number of particles per site in a large ball of
radius $m$. See Theorem \ref{T:lower_bd} for the full statement. Note that
in this result as in the rest of the paper, we will be taking limits as $m$
tends to $\infty$, recalling that $m$ is itself a function of $n$ when
applying these results to get Theorem \ref{T:long_time}.

Let us first present some further heuristic arguments for the lower bound
in \thmref{long_time}. Consider for the moment the case $d>2$, so that the
random walk migration process is transient. The first heuristic comes from
the first moment calculation and simple Green function estimates: label the
particles in an arbitrary way and let $S_i$ be the total number of
particles that ever coalesce with the $i$th particle. Observing that a
typical particle is at initial distance of order $k$ away from an order
$k^{d-1}$ particles, where $k$ ranges from $1$ up to a number of order $m$,
gives for a typical $i$
\[
\E(S_i) \asymp \sum_{k=1}^m k^{d-1} \frac{1}{k^{d-2}} \asymp m^2,
\]
where we use the fact that the probability that two particles ever coalesce
is proportional to the probability that their corresponding walks intersect
(visit the same site {\em at the same time}) (the constant comes from the
delayed coalescence dynamics). The fact that this probability is
approximately $k^{2-d}$ is a well-known Green function estimate. Since $N =
\sum_{i=1}^M 1/S_i$ gives the total number of clusters that survive forever
(with $M$ being the initial number of particles, of order $m^d$), and since
$\E(1/S_i)\geq 1/\E(S_i)$ we arrive at
\begin{equation}\label{E:expectlb}
  \E (N) \geq c m^{d-2}.
\end{equation}
While Jensen's inequality may seem crude, this does give the correct
exponents because the distribution of $S_i$ is sufficiently concentrated.
The next section contains results confirming this heuristic.

\subsection{Technical random walk lemmas}

We begin with technical results concerning random walks. Most of these are
standard yet difficult to ``pinpoint'' in the random walk literature. Let
$(S_n, n\geq 0)$ be simple symmetric random walk on $\Z^2$, started from a
point $X_0\in B(o,2m)$ which will later be chosen in a certain random
fashion (very roughly speaking, close to uniform) and recall that $m\to
\infty$. Let $(X_t, t\geq 0)$ be a continuous time random walk on $\Z^2$
obtained as $X_t:=S_{N^*_t}$ where $(N^*_t,t\geq 0)$ is a Poisson process
with rate $1$, independent of $X$. Let $S'$ be a lazy version of $S$, with
$S_0'=S_0$ and step distribution given by $\P(S_{n+1}'-S_n'= 0)=1/2$ and
$\P(S_{n+1}' -S_n'= \pm {\bf e}_i) = 1/8$, where ${\bf e}_1,{\bf e}_2$ are
the coordinate vectors, then $(S_{N^*_{2t}}', t\geq 0)$ has the same law as
$X$. We write $\P_x$ for the corresponding probability measures when
$X_0=x$.

Define $\tau_x':=\inf \{n\geq 0: S_n'=x\}$, $\tau_x:=\inf \{s>0: X_s=x\}$
to be the hitting times of $x$. Similarly, let ${\tau'}^{+}_x:=\inf \{n\geq
1: S_n'=x\}$ denote the {\em positive} hitting time of $x$. We abbreviate
$\tau'=\tau_0'$, $\tau=\tau_0$ etc.

The next result is a variation of an Erd\H{o}s-Taylor formula \cite{erdtay}
(see also \cite{coxgri}, p.~354). We assume that as $m\to \infty$,
\begin{align}
  \label{d2 hit prob hypot}
  \E\left( \frac{1}{\|X_0\|^2_+} \right) &= O \left(\frac{1}{\log
      m}\right), &
  \E\left(\log\frac{m}{\|X_0\|_+}\right) &= O(1),
\end{align}
where for any $y\in \Z$ we abbreviate $\|y\|_+:= \|y\| \vee 1$.

\begin{lemma}\label{L:d2 hit prob}
  Assume $d=2$ and fix $s>4$, assume a random $\|X_0\| \leq 2m$ satisfies
  \eqref{d2 hit prob hypot}. Then
  \begin{equation}
    \label{d2 hit prob conc}
    \P(\tau<sm^2) \asymp \frac{\log s}{\log m + \log s},
  \end{equation}
  where the constants implicit in the $\asymp$ notation depend only on
  those implicit in \eqref{d2 hit prob hypot} (and not on $s$ or $m$).
\end{lemma}

It is easy to check that $X_0$ drawn from a uniform on $B(o,2m)$ or from a
difference of two independent uniforms on $B(o,m)$ will satisfy the
hypotheses of \lemref{d2 hit prob} and therefore (\ref{d2 hit prob conc})
with universal constants (not depending on $m$) for any fixed $s>4$. Also
note that if $\P(X_0\notin B(o,2m))= 1$, under no further restriction on
the distribution of $X_0$, the upper bound on the probabilities
$\P(\tau<sm^2)$ holds with the same constant as in \lemref{d2 hit prob}.
Indeed, in order for $\tau$ to happen, the walk needs first to enter
$B(o,2m)$ at a location close to its boundary, for which the argument below
gives the required estimate.

\begin{proof}
  We estimate the above probability for any given $x \in B(o,2m)$, and then
  integrate over the law of $X_0$. Without loss of generality, assume that
  $s m^2$ is an integer. Use the ``last-exit decomposition'':
  \[
  \indic{\tau' < s m^2} = \sum_{k=1}^{s m^2-1} \indic{S_k'=0}
  \prod_{j=k+1}^{sm^2-1} \indic{S_j'\neq 0}
  \]
  together with the Markov property, to obtain
  \[
  \P_x(\tau' < s m^2) = \sum_{k=1}^{s m^2-1} \P_x(S_k'=0) \P_0(\tau'^{+} > s
  m^2 - k-1).
  \]
  We now apply a local central limit theorem and an estimate on the
  distribution of excursion length \cite{spitzer}, statement E1 on p.~167,
  (or \cite[Prop.~4.2.4]{LL-book}) and \cite{spitzer} statement P10 on
  p.~79 (or \cite[Theorem~2.1.1]{LL-book}). We find that for some universal
  sequence $e_n \xrightarrow[n\to\infty]{} 0$
  \begin{align*}
    \P_x(\tau' < s m^2) &= \sum_{k=1}^{s m^2-1}
    \frac{1}{k}e^{-\frac{2\|x\|^2}{k}} \frac{1+e_{s m^2 -k+1}}{\log(s m^2
      -k +1)} + O\left(\sum_{k=1}^{s m^2-1} \frac{1}{k \|x\|_+^2}
      \frac{1}{\log(s
        m^2 -k +1)}\right)\\
    &= \sum_{k=1}^{s m^2-1} \frac{1}{k}e^{-\frac{2\|x\|^2}{k}} \frac{1+e_{s
        m^2 -k+1}}{\log(s m^2 -k +1)} + O\left(\frac{1}{\|x\|_+^2}\right).
  \end{align*}
  Split this sum in three: For $k\leq\|x\|^2$ use $e^{-x}<x^{-1}$ to get a
  total contribution of $O(1/\log(s m^2))$. For $k > s m^2 - \sqrt{s m^2}$
  each term is at most $C/k$ so the total contribution is
  $O(1/(\sqrt{s}\cdot m))$. Finally, for the intermediate $k$'s each term
  is $\asymp 1/k\log{s m^2}$, so the total contribution is $\asymp \log(s
  m^2/\|x\|^2)/\log(s m^2)$. Thus
  \[
  \P_x(\tau' < s m^2) \asymp \frac{\log s + 2\log(m/\|x\|_+) + O(1)}{\log(s
    m^2)} + O\left(\frac1{\|x\|_+^2}\right),
  \]
  uniformly over $x \in B(o, s m/2)$. Taking expectation with respect to
  $X_0$ while using \eqref{d2 hit prob hypot} and $s>2$, yields $\P(\tau' <
  s m^2) \asymp \frac{\log s}{\log(s m^2)}$ as $m\to \infty$.

  Going back to the continuous time random walk, we have $\P( |N^*_{2t} -
  2t| > \eps t) \leq e^{-c(\eps) t}, t \geq 0$ for some $c(\eps)>0$,
  accounting for an additional error of $O(e^{-c(\eps) s m^2}) =
  o(1/\log(sm))$ in the corresponding estimate for $\tau$.
\end{proof}

We will also need later a simpler result which goes along the same lines.

\begin{lemma} \label{L:d2hit} Assume $d=2$ and $X_0=x$ is such that
  $\|x\|=m$. For all $c_1>0$, there exists $c_2>0$ which depends only on
  $c_1$ such that $\P_x( \tau_0 < c_1 m^2) \ge c_2/\log m$.
\end{lemma}

\begin{proof}
  First we note that by easy large deviations on Poisson random variables,
  it suffices to prove the same inequality with $\tau_0$ replaced by the
  discrete time $\tau'_0$. By the strong Markov property, note that if
  $K(t)$ counts the number of hits of 0 by time $t$, then for all $c>0$,
  \begin{equation}\label{d2hit1}
    \E_x\big(K(c m^2)) \le \P_x(\tau'_0 \le c m^2) \E_0(K(c m^2)).
  \end{equation}
  By the local central limit theorem,
  \begin{equation}\label{d2hit2}
    \E_0(K(c m^2))  = \sum_{k=0}^{c m^2} \P_0(S'_k = 0) \asymp
    \sum_{k=1}^{c m^2} \frac1k \sim 2c \log m.
  \end{equation}
  Also,
  \begin{equation}\label{d2hit3}
    \E_x(K(c m^2)) = \sum_{k=0}^{c m^2} \P_0(S'_k =x)
    \asymp \sum_{k=1}^{c m^2} \frac{e^{-c'\|x\|^2/(2k)}}{k}
    \ge \sum_{k=c m^2/2}^{c m^2} \frac{e^{-c'm^2/(2k)}}{k} \ge c''.
  \end{equation}
  Combining \eqref{d2hit1}--\eqref{d2hit3}, we complete the proof.
\end{proof}

If $d\geq 3$, we denote by $G_X$ the Green function of a $d$-dimensional
walk $X$. It is well-known (see e.g.~\cite{spitzer}) that
\begin{equation}
  \label{Green asymptot}
  G_X(x) \sim c\|x\|^{2-d}, \mbox{ as } \|x\| \to \infty,
\end{equation}
for some constant $c$ that depends on $d$ (here $\|x\|$ denotes the
Euclidean norm in $\Z^d$).

Let $X(\cdot),Y(\cdot)$ be independent continuous random walks in $\Z^d$
$d\geq 2$ with jump rate $1$ and with starting points uniform in $B(o,m)$.
Denote by $\sigma_t \{ X\}:=\sigma\{X(s), 0 \leq s \leq t \}$ the natural
filtration of $X$ and let $\sigma \{ X\}:= \sigma_\infty \{ X\}= \sigma\{
X(s), 0 \leq s <\infty\}$. Define the stopping time $\tau:= \inf\{ t :
X(t)=Y(t)\}$. Define the collision events by
\begin{align*}
  H_s\equiv H_s(X,Y) & := \{\tau \leq s\}, & H \equiv H(X,Y) := \{\tau <
  \infty\}.
\end{align*}

\begin{lemma}\label{L:RW_hitting_bounds}
  Let $X,Y$ be independent continuous time random walks in $\Z^d$ starting
  at uniform points at $B(o,m)$. For any $d>2$ we have
  \begin{align}
    \label{lemcon higher}
    \P(H) &\asymp m^{2-d}, \\
    \label{E var higher}
    \Var\left(\P(H|\sigma\{X\})\right) &\leq C m^{2(2-d)},
  \end{align}
  while if $d=2$, for any $t>4$ we have
  \begin{align}
    \label{lemcon}
    \P \big( H_{tm^2} \big) &\asymp \frac{\log t}{\log m+\log t}, \\
    \label{E var two}
    \Var\left(\P(H_{tm^2}|\sigma\{X\})\right) &\leq C\left(\frac{\log
        t}{\log m +\log t}\right)^2,
  \end{align}
  where $C$ and the constants in $\asymp$ relation depend only on $d$.
\end{lemma}

\begin{proof}
  Assume first that $d\geq3$. Note that the difference $X(t)-Y(t)$ is also
  a continuous time simple random walk (with a doubled rate of jumps), and
  abbreviate $G_{X-Y}=G$. It is well-known and easy to check that
  \begin{equation}
    \label{proba H Green}
  \P(H|X_0=x_0, Y_0=y_0)= \frac{G(x_0-y_0)-\indic{x_0=y_0}}{G(0)}.
\end{equation}
  Since $x\mapsto \|x\|^{2-d}$ is integrable near 0 as a function on
  $\R^d$, then \eqref{Green asymptot} implies that
  \begin{equation}\label{eq:PH_small}
    \frac{1}{\Vol B(o,m)} \sum_{y\in B(o,m)} G(x-y) \leq C m^{2-d},
    \quad \text{ for any $x\in\Z^d$}.
  \end{equation}
  If $x\in B(o,m)$, then a corresponding lower bound holds since a positive
  fraction of points in $B(o,m)$ is at distance order $m$ from $x$. Hence,
  for any $x\in B(o,m)$,
  \[
  \frac{1}{\Vol B(o,m)} \sum_{y\in B(o,m)} G(x-y) \asymp m^{2-d},
  \]
  where the constants implicit in $\asymp$ depend only on $d$. Due to
  \eqref{proba H Green}, averaging over $x\in B(o,m)$ gives that
  $\P(H)\asymp m^{2-d}$ as claimed. (It is not hard to show similarly that
  $\P(H)\sim c m^{2-d}$ for some $c$.)

  In order to show \eqref{E var higher}, introduce a third random walk $Y'$
  independent from, and identically distributed as, $X$ and $Y$. In analogy
  to $H$ define $H'=\{\exists t, X(t)=Y'(t)\}$. Given $\sigma\{X\}$, the
  events $H,H'$ are independent and have the same probability. Thus
  \begin{align*}
    \Var \P(H|\sigma\{X\})
    &\leq \E\big[ \P(H|\sigma\{X\})^2 \big] \\
    &= \E\big[ \P(H|\sigma\{X\}) \P(H'|\sigma\{X\}) \big]  \\
    &= \E\big[\P(\exists t,s: X(t)=Y(t), X(s)=Y'(s) | \sigma\{X\}) \big] \\
    &\leq 2\P(\exists t,s: t\leq s, X(t)=Y(t), X(s)=Y'(s)),
  \end{align*}
  where for the last inequality we use the symmetry between $Y$ and $Y'$.
  Denote by $\FF_\tau$ the standard $\sigma$-field generated by processes
  $X$ and $Y$ up to time $\tau$. On the event $\{\tau<\infty\}$, due to the
  strong Markov property and (\ref{proba H Green}),
  \[
  \P \big(\exists s\geq\tau: X(s)=Y'(s) ~|~ \FF_\tau \big) \leq c\E\big[
  G(X(\tau)-Y'(\tau)) ~|~ \FF_\tau \big].
  \]
  Let $Z=X(\tau)-(Y'(\tau)-Y'(0))$. Noting that $Y'(0)$ is independent from
  both $\FF_\tau$ and $Z$, we have $\E (G(Z-Y'(0)) | \FF_\tau, Z)\leq C
  m^{2-d}$, almost surely, and therefore
  \[
  \E \big(G(X(\tau)-Y'(\tau)) | \FF_\tau \big) = \E[\E \big(G(Z-Y'(0)) |
  \FF_\tau, Z\big)| \FF_\tau] \leq C m^{2-d}.
  \]
  In view of the discussion above this yields a uniform bound on $\Var
  \P(H|\sigma\{X\})$.

  If $d=2$, we proceed similarly, with $H$ replaced by $H_{tm^2}$. In
  particular, \lemref{d2 hit prob} gives the asymptotics of $\P(H_{tm^2})$.
  For the conditional variance estimate, one obtains as above
  \[
  \Var \P(H_{tm^2}(X,Y)|\sigma\{X\}) \leq 2\P[\indic{\tau <t}
  \P(H_{tm^2}(X'',Y'')|\FF_\tau)],
  \]
  where $X'',Y''$ are independent random walks started from $X(\tau)$ and
  $Y'(\tau)$, respectively, and otherwise independent of $\FF_\tau$. The
  result follows as before, since by \lemref{d2 hit prob},
  $\P(H_{tm^2}(X'',Y'')|\FF_\tau) \asymp \frac{\log{t}}{\log{m}+\log{t}}$.
\end{proof}

\subsection{Proof of the lower bound}
\label{SS:lowerbound2}

We return to the spatial coalescent. Let $\Lambda$ be an arbitrary finite
measure on $(0,1)$. Consider a spatial coalescent with initial
configuration $X(0)$ that stochastically dominates i.i.d.\ Bernoulli random
variables with mean $p>0$ in $B(o,m)$ (we make no assumptions on the
initial configuration outside of $B(o,m)$). With a slight abuse of
notation, we write $N^m(t)$ in this section for the total number of
particles at time $t$, and we define $N \equiv N^m =\lim_{t\to\infty}
N^m(t)$ be the number of particles that survive to time $\infty$.

\begin{thm}\label{T:lower_bd}
  Consider the spatial coalescent with initial state dominating Bernoulli
  variables in $B(o,m)$. If $d>2$, then there exist a constant $a>0$ such
  that
  \[
  \P(N>am^{d-2}) \xrightarrow[m\to\infty]{} 1.
  \]
  If $d = 2$, then there exists a constant $a>0$ such that, for any $t>4$,
  \[
  \P \left( N^m(tm^2) > a \, \frac{\log m}{\log t} \right)
  \xrightarrow[m\to\infty]{} 1.
  \]
\end{thm}

Note that, since the total number of particles is non-increasing, the lower
bound in the $d=2$ case holds for any $t>1$ with modified constant $a$ (or
with $\log(2+t)$ in place of $\log t$ for any positive $t$).

We begin with a lemma stating a similar result for a simpler initial
condition and with an ``instantaneous'' coalescent mechanism, where two
particles coalesce as soon as they visit the same site. This model is
called {\em coalescing random walks} (CRW). Afterwards we couple the two
models to obtain \thmref{lower_bd}.

\begin{lemma}\label{L:lower_bd_iid}
  Consider a system of $s$ coalescing random walks, such that their initial
  positions are i.i.d.\ uniform points in $B(o,m)$, where
  \[
  s \equiv s(a) = \begin{cases}
    a m^{d-2},& d\geq 3,\\
    a \log m, & d=2.
  \end{cases}
  \]
  Let $Z(t)$ denote the total number of particles at time $t$ and let $Z=
  \lim_{t \to \infty} Z(t)$. If $d>2$, then for some $a>0$ we have
  $\P(Z>am^{d-2}/4) \xrightarrow[m\to\infty]{} 1$. \\
  If $d = 2$, then for some $a$ and all $t>4$, we have $\P \left( Z(tm^2) >
    a \frac{\log m+\log t}{4 \log t} \right) \xrightarrow[m\to\infty]{} 1$.
\end{lemma}

\begin{proof}
  We use the following explicit construction of the CRW model with the
  given initial condition: Let $(X_i(t), t\geq 0), i=0,1,\dots,s-1$ be a
  family of i.i.d.\ (non-coalescing) random walks, such that for each $i$,
  $X_i(0)$ is uniform in $B(o,m)$. At time $0$, each block contains a
  single particle that is assigned a unique label in $\{0,1,\dots, s-1\}$.
  While present in the system, the particle (or block of particles)
  carrying label $i$ follows the trajectory of $X_i$. If the trajectories
  of blocks labeled $i$ and $j$ ever intersect, they instantaneously merge
  into a new block that inherits the smaller label $i\wedge j$.

  Consider first the case $d>2$. For each pair $i,j$ let $A_{i,j} :=
  \{\forall u\geq 0: X_i(u) \neq X_j(u)\}= H(X_i,X_j)^c$. Then on $A_{i,j}$
  the blocks carrying labels $i$ and $j$ cannot merge as a consequence of a
  single coalescence event, but might merge due to a collection of
  coalescence events involving lower indexed particles. However, on the
  event
  \begin{equation}\label{EAk}
    A_k := \bigcap_{i<k} A_{k,i},
  \end{equation}
  the block carrying label $k$ stays in the system indefinitely.

  Consider the filtration $\FF_k=\sigma\big\{X_i(\cdot), i\leq k\big\}$.
  Define $p_k = \P(A_k|\FF_{k-1})$, and note that $p_0=1$. The random
  variables $\{p_k\}$, are a non-increasing sequence of random variables.
  To see this we use the fact that the random walks are independent and so
  \[
  p_k \leq \P\left(\bigcap_{i=0}^{k-2} A_{k,i}|\FF_{k-1}\right) =
  \P\left(\bigcap_{i=0}^{k-2} A_{k-1,i}|\FF_{k-2}\right) = p_{k-1}, \mbox{
    almost surely}.
  \]

  Next, define events
  \[
  B_k = A_k \cup \{p_k<1/2\},
  \]
  and note that
  \[
  \P(B_k|\FF_{k-1}) = \begin{cases}
    1, & p_k<1/2, \\
    p_k, & p_k\geq 1/2.
  \end{cases}
  \]
  Consider the martingale
  \[
  M_k = \sum_{i=0}^k \indica{B_i} - \P(B_i|\FF_{i-1}).
  \]
  Note that $M_k$ has increments with variance bounded (crudely) by 1. Thus
  $\Var M_s<s$ (here $s$ is the initial total number of blocks) and, by
  Markov's inequality,
  \[
  \P(|M_s| > s/4)\leq \frac{s}{(s/4)^2} = \frac{16}{s}.
  \]
  However, $\P(B_k|\FF_{k-1})\geq1/2$, so by the definition of $M$, we find
  \begin{equation}\label{EsumBlowerbd}
    \P\left(\sum_{i<s} \indica{B_i} < s/4 \right)\leq \frac{16}{s}
    \xrightarrow[m\to\infty]{} 0.
  \end{equation}

  Since $p_k$ is non-increasing and since on the event $\{p_k\geq1/2\}$ the
  events $A_k$ and $B_k$ coincide, we realize that on the event $\{p_s\geq
  1/2\}$
  \[
  \sum_{i<s} \indica{A_i} = \sum_{i<s} \indica{B_i}.
  \]
  Thus if we prove that
  \begin{equation}\label{Epclaim}
    \P(p_s<1/2) \xrightarrow[m\to\infty]{} 0,
  \end{equation}
  then \eqref{EsumBlowerbd} would imply the lemma. To this end we show that
  $p_s$ is bounded below by a random quantity that is concentrated above
  $1/2$, via second moment estimates. Specifically, from the definition
  \eqref{EAk} we have
  \[
  1 - p_s \leq \sum_{i<s} \P(A_{s,i}^c| \FF_{s-1}) = \sum_{i<s}
  \P(A_{s,i}^c|\sigma\{X_i\}),
  \]
  where the last identity is due to independence of $\sigma\{X_i\}$ for
  different $i$'s. Moreover, $\{\P(A_{s,i}^c|\sigma\{X_i\}),
  i=0,\ldots,s-1\}$ is an i.i.d.\ family of random variables. Using
  \eqref{lemcon higher},
  \[
  \E\left(\sum_{i<s} \P(A_{s,i}^c|\sigma\{X_i\}) \right) < s \cdot C
  m^{2-d} \leq C a.
  \]
  We choose $a=1/(4C)$ so that this expectation is at most $1/4$. Due to
  \eqref{E var higher},
  \[
  \Var\left(\sum_{i<s} \P(A_{s,i}^c|\sigma\{X_i\}) \right) \leq s\cdot C
  m^{2(2-d)} \to 0.
  \]
  so the sum is concentrated near its mean, and \eqref{Epclaim} follows.

  \medskip

  In the case $d=2$, the proof is almost identical. We take $s = a \log m$
  and $a < 1/(4C \log t)$, where $C$ is the constant that appears in
  \eqref{lemcon}. The event $A_{i,j}$ is accordingly redefined as $A_{i,j}
  := H_{tm^2}(X_i,X_j)^c$. Otherwise, the argument proceeds exactly as
  above, with \eqref{lemcon}, \eqref{E var two} used in place of
  \eqref{lemcon higher}, \eqref{E var higher}.
\end{proof}

\begin{proof}[Proof of \thmref{lower_bd}]
  The idea is to couple the spatial coalescent $X$ with a system of
  coalescing random walks, denoted $X\crw$, with an initial state of $s$
  particles at i.i.d.\ sites, uniform in $B(o,m)$. We first argue that it
  is possible to couple the initial states so that w.h.p.\ $X\crw_v(0)\leq
  X_v(0)$ (at every vertex). Indeed, in $X\crw$, there are $N\crw(0)\leq s$
  occupied sites (since there may be repetitions) and given $S$, these
  sites are uniformly sampled from the ball $B(0,m)$ without replacement.
  On the other hand, $X(0)$ dominates a Bernoulli configuration on
  $B(o,m)$, hence $X(0)$ has at least $\Bin(\#B(o,m),p)$ particles sampled
  without replacement. Since $\P(\Bin(\#B(o,m),p) > s) \to 1$, this holds.

  The second step of the proof is that if the initial configurations
  satisfy $X\crw_v(0) \leq X_v(0)$ for all $v$, then there is a coupling of
  the processes so that
  \begin{equation}\label{couplingCRW}
    X_v(t) \ge X\crw_v(t), \ \ t\ge 0, v \in V.
  \end{equation}
  To see this, observe that by the consistency property of spatial
  $\Lambda$-coalescent it suffices to prove the result assuming that
  $X_v(0) = X\crw_v(0)$ for all $v \in V$. In this case,
  \eqref{couplingCRW} follows easily by induction on the number of
  particles: Just apply the consistency property of spatial
  $\Lambda$-coalescents, after the first time that two particles occupy the
  same site. (This idea is further exploited in Lemma \ref{L:sumdom2}.)

  Finally, \thmref{lower_bd} follows by \lemref{lower_bd_iid}.
\end{proof}

\subsection{Concentration of the number of particles}
\label{SS:concentration}

The main result of this section is a concentration result for the number of
particles alive at a certain time. This provides a soft alternate route for
the lower-bound on the long-time behavior of the spatial coalescent, as we
briefly explain.

\begin{thm}\label{T:lt-var}
  Fix $t>0$, and consider a spatial Kingman coalescent started from some
  arbitrary configuration containing a finite number of particles. Then we
  have
  \[
  \Var (N(t))\leq \E N(t) .
  \]
\end{thm}

\begin{proof}
  The tool used here again is a comparison to the coalescing random walk
  model, where particles coalesce immediately upon meeting. We denote by
  $(X\crw(t),t\ge 0)$ a system of instantaneously coalescing random walks
  started from a certain set of vertices $A$ in a graph $G=(V,E)$, to be
  chosen suitably later, and let $N\crw(t)$ denote the total number of
  particles at time $t$. The proof is based on Arratia's correlation
  inequality \cite[Lemma~1]{arratia}, which states that
  \begin{equation}\label{Ncorr-A}
    \E X\crw_x(t) X\crw_y(t)  \leq  \E X\crw_x(t) \cdot \E X\crw_y(t).
  \end{equation}
  Thus at any time, any two sites are negatively correlated. This
  inequality holds not just for the process on $\Z^d$, but on any edge
  weighted graph.

  We now remark that the spatial Kingman coalescent on $\Z^d$ can be
  approximated by a system of instantly coalescing random walks on a larger
  graph. For any integer $N$ such that $N>n$ (the initial number of
  particles), consider the graph $G_N=(V,E)$ with vertices
  $V=\Z^d\times\{1,\dots,N\}$. The edges of $G_N$ are of two types. If
  $x\sim y$ in $\Z^d$ then there is an edge between $(x,i)$ and $(y,j)$
  with weight $\rho/N$. Additionally, there is an edge with weight $1/2$
  between $(x,i)$ and $(x,j)$ for any $x,i,j$. Call the set
  $x\times\{1,\dots,N\}$ a cluster. Clusters correspond to vertices of
  $\Z^d$ in a natural way. The $\Z^d$ coordinate of a continuous time
  random walk on $G_N$ is a continuous time random walk on $\Z^d$ with jump
  rate $\rho$. However, two walks may be present in the same cluster and
  not meet. It is clear that as long as two random walks are in the same
  cluster they will meet at rate one (since each may jump into the vertex
  occupied by the other).

  The probability of two random walks meeting when one jumps from one
  cluster to another is of order $1/N$. Thus as long as the number of
  particles is negligible compared to $N$, the projection onto $\Z^d$ of
  the coalescing random walks $X\crw_N$ on $G_N$ is close to the spatial
  Kingman coalescent on $\Z^d$. As $N\to\infty$, the projection of
  $X\crw_{N}(t)$ converges to $X(t)$ (in the sense of vague convergence,
  identifying $X_v$ and the projection of $X\crw_N$ to point measures on
  $\Z^d$).
  More precisely, for an initial configuration $X(0)$ of particles on
  $\Z^d$, we define a set $A\subset V_N$ by choosing for each $v$
  (arbitrarily) $X_v(0)$ particles from the cluster of $v$. Let
  $X\crw_{N}(t)$ be the process of coalescing random walks on $G_N$ started
  with this configuration. Then if $M_N\crw(t)$ denote the total number of
  particles of $X\crw_N(t)$,
  \begin{align*}
    \E(M_N\crw(t)^2) &= \sum_{x\in V_N} \E X\crw_{N,t}(x)
    + \sum_{x\neq y\in V_N} \E X\crw_{N,t}(x) X\crw_{N,t}(y) \\
    &\leq \E M_N\crw(t)
    + \sum_{x\neq y\in V_N} \E X\crw_{N,t}(x) \E X\crw_{N,t}(y) \\
    &\leq \E M\crw_N(t) + (\E M\crw_N(t))^2.
  \end{align*}
  Thus for any $N$ we have $\Var M\crw_N(t) \leq \E M\crw_N(t)$. By
  dominated convergence (since all processes have at most $n$ particles) we
  see that
  \begin{align*}
    \lim_{N\to\infty} \E M\crw_N(t) &= \E N(t) & \lim_{N\to\infty} \E
    M\crw_N(t)^2 &= \E N(t)^2,
  \end{align*}
  and the result follows.
\end{proof}

As a simple corollary of this result, we obtain an alternate proof of
Theorem~\ref{T:lower_bd}. We have already seen in (\ref{E:expectlb}) that
$\E(N(\infty)) \ge c m^{d-2}$ for some $c>0$ if $d \ge 3$ (this argument is
a simple Green function estimate, and is easy to adapt to the case $d=2$).
Applying \thmref{lt-var} concludes the proof.

\medskip

It would be also possible to derive a lower-bound on the expected number of
particles in a system of instantaneously coalescing random walks at time
$tm^2$, starting from a set $A$ which dominates i.i.d. Bernoulli random
variables with mean $p>0$, using technology from coalescing random walks.
We briefly outline the steps needed to do this. First, starting from a
configuration where there is a particle at every site of $\Z^d$, and using
a famous result of Bramson and Griffeath \cite{bg80} on the asymptotic
density of particles, we conclude that about $c m^{d-2}$ such particles are
in a region of volume $C m^d$ for some large $C>0$ to be chosen suitably.
If we treat the particles that started outside of $A$ as ghosts, we are
then led to estimate the number of ghost particles among those $cm^{d-2}$.
For this, one can use the duality with the voter model (see \cite{liggett})
and \cite[Lemma~4]{merle}, which gives good control on the probability that
the voter model escapes a ball of radius $\sqrt{t}$, for large $t$.

\section{Upper bound for the number of survivors}
\label{SS:upper bound}

Assume that $\Lambda$ is a finite measure on $[0,1]$ such that for some
$a_0>0$, we have
\begin{equation}\label{condUB}
  \lambda_n \ge a_0 n \ \ \text{ for all } n\ge 2,
\end{equation}
where $\lambda_n = \sum_{k=2}^n \lambda_{b,k}$ is the total merger rate
when there are $n$ particles. Note that most coalescents which come down
from infinity satisfy \eqref{condUB}, in particular, if $\Lambda=
\delta_{\{0\}}$ (the Kingman case) then \eqref{condUB} holds since
$\lambda_n={n \choose 2}$, and if $\Lambda$ has the regular variation
property of \eqref{reg-var}, then \eqref{condUB} holds by
\lemref{lambda_asymp}.

Our goal here is to prove the following result.

\begin{thm}\label{T:upper_bd}
  Fix $C_0\in (0,\infty)$ and $\delta>0$, and consider the spatial
  $\Lambda$-coalescent in $\Z^d$ satisfying (\ref{condUB}), started from a
  configuration of at most $C_0m^d$ particles located in $B(o,m)$, and no
  particles in $\Z^d\setminus B(o,m)$. There exists $C=C(\delta,C_0)$, such
  that if $d>2$ then
  \[
  \P\left( \tN(\delta m^2) < Cm^{d-2} \right) \xrightarrow[m\to\infty]{} 1,
  \]
  while, if $d=2$,
  \[
  \P\left( \tN(\delta m^2) < C \ln m \right) \xrightarrow[m\to\infty]{} 1.
  \]
\end{thm}

Note that when $d>2$ this order of magnitude bound is sharp, since
\thmref{lower_bd} shows $\tN(\infty)\geq cm^{d-2}$. For $d=2$, due to
recurrence, $\tN(\infty)=1$, almost surely.

The idea behind the proof is a comparison of the spatial system to a mean
field approximation. The actual argument is based on a somewhat technical
construction so we start with a non-technical overview. Recall the
comparison with ODE described in (\ref{rho_heur}): if at time $t$ the
density of particles averaged over some ball is $\rho(t)$ (typically
small), then we approximate the spatial coalescent with the mean-field
model where the coalescence rate per particle is $\rho(t)$ at time $t$,
leading to the differential equation
\[
\frac{d}{dt} \rho(t) = -\rho^2(t)/2, \ \ t\geq s.
\]
Hence $\rho(t)^{-1} = c+(t-s)/2$ and therefore $\rho(t) = \frac{2}{t - s +
  2\rho(s)^{-1}}$, $t\geq s$. Provided that all the particles in the
spatial coalescent configuration are located in the ball of radius $m$
during the whole interval $[s,t]$ (and that the above approximation is
valid) then their total number is approximately $C m^d \rho(t)$. In turn,
this approximation remains valid as long as the particles remain inside a
ball centered at the origin with radius of order $m$, i.e.\ up to time of
order $m^2$. At times of order $m^2$, the number of remaining particles is
of order $m^{d-2}$.

A key difficulty of the approach outlined above comes from the fact that
some particles diffuse away from the densest regions relatively early in
the evolution, which might enable them to survive longer. To account for
such ``runaways'', we adopt a {\em multi-scale approach}, bounding at each
stage the number of particles that ``escape''. This is done in
\lemref{ub_induction}. \lemref{low_density} provides the estimates on the
number of non-escaping particles at each stage.

To justify the comparison of the spatial process with the mean field
process we average over small time intervals (cf.~\lemref{time_average}
below). This is necessary since at any given time it is possible that no
vertex contains more than a single particle, in which case the immediate
rate of coalescence is 0. However, the system is unlikely to stay in such
states long enough to hinder the approximation. Indeed,
\lemref{time_average} implies that the average rate of coalescence is (up
to constants) as predicted by the mean field heuristic. The multiplicative
constants are inherent to the spatial structure, and it seems difficult to
compute them.

\subsection{Preparatory lemmas}

Our first step is a comparison lemma between the spatial
$\Lambda$-coalescent $X$ and a slower spatial coalescent. We then consider
a possibly more general spatial coalescent process $\{(\bar X_v(t),t\ge
0)\}_{v \in V}$. If the process consists initially of $n$ particles labeled
by $[n]=\{1,\ldots, n\}$, a configuration consists as usual of labeled
partitions of $[n]$, where the label of a block corresponds to its location
on $V$. Equivalently, a configuration $\bar x=(\bar x_v)_{v \in V}$ may be
thought of as giving the list of blocks (referred to as particles) present
at each particular site $v \in V$. We will also sometimes abuse notation
and denote by $X_v(t)$ the number of particles (i.e., blocks) present at
time $t$ and at position $v$. We assume that particles perform independent
continuous-time simple random walks with jump rate $\rho$, and that there
exists a family of real numbers $\bar \lambda_{\bar x, S}$ such that for
all configuration $\bar x = (\bar x_v)_{v \in V}$, all $v \in V$, any
particular subset $S$ of all blocks present at $v \in V$ coalesces at an
instantaneous rate $\bar \lambda_{\bar x, S}$, if the current configuration
is $\bar x$. Moreover, coalescence events at different sites occur
independently of one another, and are independent of the migration. We now
make the following assumption on the family of rates $\bar \lambda_{\bar x
  ,S}$: if $v \in V$ and $\bar x_v$ contains $n \ge 2$ particles, then for
every $2\le k \le n$, we have:
\begin{equation}\label{eq:dom_rates}
  \sum_{S: |S|\geq k} \bar\lambda_{\bar x,S} \leq \sum_{\ell \ge k}
  \binom{n}{\ell} \lambda_{n,\ell},
\end{equation}
where $\lambda_{n,k}$ is the coalescence rate of any particular subset of
size $k$ in a $\Lambda$-coalescent. The idea behind \eqref{eq:dom_rates} is
that if $X$ and $\bar X$ have the same number of particles at time $t$,
then $X(t+\eps)$ is stochastically dominated by $\bar X(t+\eps)$.

\begin{lemma}\label{L:sumdom2}
  Consider a $\Lambda$-coalescent $X$ and a coalescent process $\bar X$
  such that (\ref{eq:dom_rates}) holds, and $X_v(0)\leq \bar X_v(0)$ for
  all $v$. Then there is a coupling of the processes $X$ and $\bar X$ such
  that $X_v(t)\leq \bar X_v(t)$ holds for all $v \in V$ and $t \ge 0$.
\end{lemma}

\begin{proof}
  By the consistency of spatial $\Lambda$-coalescents, it suffices to prove
  the result when $X_v(0) = \bar X_v(0)$ for all $v \in V$. We associate
  each particle of $X$ with a particle of $\bar X$ and let them perform the
  same random walks as long as there are no coalescence events. A
  consequence of \eqref{eq:dom_rates} is that it is possible to couple the
  processes so that if $X_v=\bar X_v$ then the coalescence events of $X$
  dominate those of $\bar X$, that is, any coalescence event in $\bar X$
  occurs at the same time as an event in $X$ involving at least as many
  particles.

  The proof now proceeds by induction on the total number of particles,
  which are allowed to be distributed arbitrarily. By the above remark, we
  may couple the processes $X$ and $\bar X$ so that the domination holds up
  to and including the first time $t_0$ of a coalescence event, which
  involves particles from $X$ and possibly from $\bar X$. Assume that $\bar
  X$ also experiences a coagulation event at this time. (Else, we can
  artificially retain particles in $X$ that were supposed to coagulate at
  time $t_0$. By the consistency property, this may only increase the
  process $X$ stochastically.)

  We now use the induction hypotheses to construct processes $(X'(t),t\ge
  t_0)$ and $(\bar X'(t), t\ge t_0)$ with initial configuration $X'(t_0) =
  \bar X'(t_0) = \bar X(t_0)$ such that $X'_u(t) \le \bar X'_u(t)$ for all
  $t \ge t_0$. We can define $\bar X(t)=\bar X'(t)$ for $t>t_0$, and by
  consistency of the spatial $\Lambda$-coalescents, we extend the coupling
  to $X$ for $t>t_0$ so that $X_u(t) \le X'_u(t)$ for all $u\in V$, which
  proves the claim.
\end{proof}

\begin{rem*}
  This lemma holds for more general spatial coalescents: e.g., the
  instantaneous coalescence rates $\lambda_{\bar x, S}$ could be allowed to
  be arbitrary path-dependent (i.e., $\mathcal{F}_t$-measurable at time
  $t$), almost surely nonnegative and finite random variables. The only
  crucial assumption is that (\ref{eq:dom_rates}) holds uniformly.
\end{rem*}

We now apply Lemma \ref{L:sumdom2} to the situation which is particularly
useful in our setting. Recall that we are considering a spatial
$\Lambda$-coalescent for which (\ref{condUB}) holds. Assume that initially
there are $N$ particles, and let $\{X_v(t),t\ge 0\}_{v \in V}$ denote the
number of particles of this process as a function of time and space.

Let $\pi$ be a partition of $\{1, \ldots, N\}$. We refer to the blocks of
$\pi$ as {\em classes}. Let $\{\bar X_v(t),t\ge 0\}_{v \in V}$ denote a
process where classes evolve independently of one another, and particles
within each class evolve according to a spatial $(\bar
\Lambda)$-coalescent, where $\bar \Lambda$ will be specified soon. That is,
particles move as continuous-time simple random walks with rate $\rho$ and
coalesce when they are on the same site and from the same class according
to a $\bar \Lambda$-coalescent.

\begin{lemma}
  \label{L:sumdom3}
  Assume that the blocks of $\pi$ are all of size 1 or 2, and that $\bar
  \Lambda = (a_0/\lambda_2) \Lambda$, where $a_0$ is the constant of
  (\ref{condUB}) and $\lambda_2=\lambda_{2,2}$ is the pairwise coalescence
  rate. Assume also that $X_v(0)\leq \bar X_v(0)$ for all $v$. Then there
  is a coupling of the processes $X$ and $\bar X$ such that $X_v(t)\leq
  \bar X_v(t)$ holds for all $v \in V$ and $t \ge 0$.
\end{lemma}

\begin{proof}
  Observe first that our process $\bar X$ is of the type described above
  \lemref{sumdom2}, so that it suffices to establish \eqref{eq:dom_rates}.
  Note however that if a configuration $\bar x$ contains $n$ particles at
  site $v$, and $S$ is a subset of particles with $|S|=k$ and $2\le k \le
  n$, we have $\bar \lambda_{\bar x, S} = 0$ for $k\ge 3$, while if $k=2$,
  $\lambda_{\bar x, S} = 0$ when the particles of $S$ are not of the same
  class, and if they are of the same class, $\lambda_{\bar x, S} =
  (a_0/\lambda_2) \lambda_2 = a_0$. Since there are at most $n$ subsets of
  particles that are allowed to coalesce, we have
  \begin{align*}
    \sum_{S:|S|\ge 2} \bar \lambda_{\bar x,S} \le n a_0 \le \lambda_n =
    \sum_{k=2}^n {n \choose k} \lambda_{n,k},
  \end{align*}
  which proves (\ref{eq:dom_rates}), and completes the proof.
\end{proof}

\begin{lemma}\label{L:time_average}
  Fix $c_0,C_0$, and consider a spatial $\Lambda$-coalescent satisfying
  (\ref{condUB}) with $\tN(0)$ particles all inside $B(o,R)$. Let $\rho(t)
  = \frac{\tN(t)}{R^d}$ be the inverse density, and assume $\rho(0) \in
  (c_0 R^{-2}, C_0)$. Denote $\tau=\rho(0)^{-2/d}$. Then for $d>2$ we have
  \[
  \P\Big(\rho(\tau)^{-1} < \rho(0)^{-1} + c_1\tau\Big) < \exp\left( -c
    R^{(d-2)^2/d} \right),
  \]
  where $c_1$ depends only on $d,c_0,C_0,a_0$.

  If $d=2$ we have
  \[
  \P\Big(\rho(\tau)^{-1} < \rho(0)^{-1} + \frac{c_1\tau}{\log\tau}\Big) <
  \exp\left( -\frac{c R^2}{\tau \log \tau} \right).
  \]
\end{lemma}

\begin{proof}
  We first argue that for some $C=C(d)$, it is possible to find at least
  $N^*(0)/4$ disjoint pairs in the set of initial particles, so that for
  each pair the initial distance between its particles is at most $C
  \rho(0)^{-1/d}$ (for large $\rho(0)$, the particles forming such a pair
  are initially located at the same site). To achieve this, cover $B(o,R)$
  with $N^*(0)/2$ (disjoint) boxes of diameter $C R N^*(0)^{-1/d} =
  C'\rho^{-1/d}$ (this is possible for some $C$). Within each box match as
  many pairs as possible in an arbitrary manner. This leaves at most one
  unpaired particle in each ball, so at least $N^*(0)/2$ are matched, with
  all distances bounded as claimed. Refer to two particles forming a pair
  as ``partners''.

  Consider the coupling from Lemma \ref{L:sumdom3}, where $\pi$ corresponds
  to the partitioned formed by identifying particles with their partners
  (which therefore contains only singletons or doubletons). Let $Z'$ be the
  total number of coalescence events in the process $\Pi'$ where
  coalescence events involving members of different classes are not allowed
  and occur at rate $a_0$ when they are. Lemma \ref{L:sumdom3} implies that
  $Z'\preceq Z$, in the sense of stochastic domination. Hence, it suffices
  to prove the claimed bounds for $Z'$. The advantage of considering $\Pi'$
  instead of $\Pi$ is that different pairs of partners evolve
  independently.

  From this point on, the arguments for the cases $d=2$ and $d>2$ differ.
  In dimensions $d>2$, by our assumptions, $\tau>c_0$ for some $c_0$. The
  probability that random walkers started at distance at most $\rho^{-1/d}$
  meet before time $\tau/2=\rho^{-2/d}/2$ is at least $c\rho^{(d-2)/d}$. On
  this event, there is probability bounded from 0 that they coalesce before
  time $\tau$. Thus the number of partners that coalesce by time $\tau$
  dominates a $\Bin(\tN(0)/4, c\rho^{(d-2)/d})$ random variable. This
  random variable has expectation $c\tN(0)\rho^{(d-2)/d} \geq c
  R^{(d-2)^2/d}$. The bound in the lemma is the probability that this
  random variable is less than half its expectation.

  Finally, if the number of coalesce events is at least
  $c\tN(0)\rho(0)^{(d-2)/d} = c\tN(0)\rho(0)\tau$ then
  \[
  \rho(t)^{-1} = \frac{R^d}{\tN(t)} \geq \frac{R^d}{\tN(0)(1-c\rho(0)\tau)}
  \geq \frac{R^d}{\tN(0)} (1 + c\rho(0)\tau) = \rho(0)^{-1} + c\tau.
  \]

  In the case $d=2$, each pair coalesces with probability at least
  $c/\log\tau$ (along the same lines) by Lemma \ref{L:d2hit}. As above, the
  number of coalesce events is at least $c R^2/\tau\log\tau$ except with
  probability $e^{-c R^2/\tau\log\tau}$. On this event, a similar
  computation gives
  \[
  \rho_t^{-1} \geq \rho_0^{-1} + \frac{c}{\tau\log\tau}. \qedhere
  \]
\end{proof}

\begin{lemma}\label{L:low_density}
  Let $S_t$ denote the number of particles (in the spatial coalescent) that
  remain in $B(o,R)$ during the whole interval $[0,t]$. In particular,
  $S_0= \sum_{v\in B(o,R)} X_v(0)$. Fix $C_0>0$, and assume that $S_0=n<C_0
  R^d$ and that $2<t<R^2$. Then for some $C$ depending only on $d$,
  \[
  \P\left( \frac{S_t}{R^d} > \frac{C}{t} \right) < c R^4 \exp\big(-c
  R^{(d-2)^2/d} \big) \qquad\text{if $d>2$},
  \]
  and
  \[
  \P\left( \frac{S_t}{R^2} > \frac{C\log t}{t} \right) < (\log t)^{-1}
  \qquad\text{if $d=2$}.
  \]
\end{lemma}

\begin{proof}
  As usual, the case $d>2$ is considered first. The previous lemma can be
  formulated as follows: The process $\{\rho_t^{-1}\}$ is unlikely to spend
  more than $u^{-2/d}$ units of time in the interval $[u,u + c_1
  u^{-2/d}]$. Note that $S_t$ can only decrease faster than $\tN(t)$, so
  this will also hold for the modified density $\rho(t)=\frac{S_t}{R^d}$.

  We apply this to the following sequence of intervals. Let $u_0=\rho(0)$
  and $u_{k+1} = u_k + c_1 u_k^{-2/d}$. Let $K$ be minimal with $u_K >
  t/c_1$. As long as $u_j<t/c_1$ the increment is at least $c t^{-2/d}$. It
  follows that $K<C t^{1+2/d} < R^4$. If the process does not spend more
  than $u_k^{-2/d}$ time in $[u_k,u_{k+1}]$ then the time before
  $\rho^{-1}$ exceeds $t/c_1$ is at most $t$. The probability that this
  fails to hold is at most $R^2 \exp\big(-c R^{(d-2)^2/d} \big)$.

  This works provided $K>1$, or equivalently $t \geq \rho(0)^{-2/d}$. If $t
  < \rho(0)^{-2/d}$ then we have
  \[
  S_t \leq S_0 = \rho(0) R^d \leq \frac{R^d}{t^{d/2}} < \frac{R^d}{t}.
  \]

  \medskip

  In the case $d=2$, we instead have $u_{k+1} = u_k + \frac{c_1 u_k}{\log
    u_k}$. It is not hard to see that
  \begin{equation}\label{uk}
    u_k \asymp e^{\sqrt{2c_1 k + c}}
  \end{equation}
  for some $c$ depending on $u_0$. To this end, note that $u$ is increasing
  and hence is dominated by the solution of the ODE $f'=c_1 f/\log f$ (at
  least once $u$ is large enough that $u/\log u$ is increasing). This ODE
  is solved by $f=e^{\sqrt{2c_1 x+x}}$, giving the upper bound on $u$. For
  the other direction, note that once $u_k$ is large $u_{k+1}/u_k$ is close
  to $1$. This implies that $u$ dominates a solution of $f'=(c_1-\eps)
  f/\log f$.

  \lemref{time_average} tells us that $\rho_t^{-1}$ is unlikely to spend
  more than $u_k$ units of time in $[u_k, u_{k+1}]$, and the probability of
  this unlikely event is at most
  \[
  p_k = e^{-c R^2/u_k\log u_k}.
  \]
  Let $K$ be such that $u_K>\alpha t/\log t$, with $\alpha$ small to be
  determined soon. Note that now the failure probability for the last
  intervals is of order 1, so a union bound does not work. However,
  \lemref{time_average} tells us more. If the process $\rho_t^{-1}$ fails
  to exceed $u_{k+1}$ in the next $u_k$ units of time, then by the Markov
  property and Lemma \ref{L:time_average} again, it gets a fresh chance to
  do so in the next $u_{k}$ units of time. Therefore, the number of
  attempts is smaller than a geometric random variable with success
  probability $1-p_k$. It follows that the total time spent in $[u_0,u_K]$
  is stochastically dominated by $Q := \sum_{k=1}^{K-1} u_k \Geom(1-p_k)$,
  where the geometric random variables are independent. The probability we
  wish to bound is therefore at most $\P(Q>t)$. By making $\alpha$ small we
  can guarantee $p_k<1/2$ for all $k$, so that the geometric variables are
  typically small.

  More precisely, from (\ref{uk}) it follows that $K\asymp \log^2 t$ and
  therefore that
  \[
  \sum_{k<K} u_k \asymp \int_0^K e^{\sqrt{2c_1 t}} dt \asymp \sqrt{K} u_K
  \asymp \alpha t,
  \]
  hence for small enough $\alpha$ we have $\E Q \leq 2\sum_{k=1}^{K-1} u_k
  \leq t/2$. Similarly, we can compute
  \[
  \Var Q \leq C \sum_{k<K} u_k^2 < \frac{C\alpha^2 t^2}{\log t}.
  \]
  The lemma now follows from Chebyshev's inequality and choosing $\alpha$
  small enough that $4\alpha^2 C \le 1$.
\end{proof}

The following is a fairly standard fact which follows easily from the
optional stopping theorem and Doob's inequality:

\begin{lemma}\label{L:gaussian}
  If $X_t$ is a continuous time random walk on $\Z^d$ then for all $x \in
  \Z^d$,
  \[
  \P\left(\sup_{s\leq t } \|X_s\|_2 \ge x \right) \leq C e^{-c x^2/t},
    \]
  where $c,C$ depend only on $d$.
\end{lemma}

\begin{lemma}\label{bin_dev}
  If $W = \Bin(n,p)$, and $\Delta > 0$, then
  \begin{equation}\label{eq:bin_dev}
    \P(W > 2 n p + \Delta) < \Delta^{-1/2}.
  \end{equation}
\end{lemma}

\begin{proof}
  If $\Delta > (n p)^2$, then Markov's inequality gives $\P(W>\Delta) \leq
  n p/\Delta < \Delta^{-1/2}$. If $\Delta \geq (n p)^2$, one can use
  Chebyshev's inequality to obtain $\P(W > 2 n p + \Delta) \leq \P(W - \E W
  > n p) \leq (1-p)/(n p) < \Delta^{-1/2}$.
\end{proof}

\subsection{Completing the proof}

\lemref{low_density} is almost sufficient to deduce \thmref{upper_bd}. The
missing piece is to account for the particles that ``escape'' from the ball
under observation. We accomplish this by partitioning the time interval
$[0,m^2]$ into several segments and applying \lemref{low_density} to each
segment. More precisely, let $K = K(m) = \lfloor \log\log m \rfloor$, and
consider the process at a particular sequence of times given by
\[
t_k = \begin{cases} 0 & k=0, \\ e^{k-K} m^2 & k=1,\dots,K.
\end{cases}
\]
Thus $t_1 \asymp m^2/\ln m$, and the sequence increases geometrically up to
$t_K=m^2$. At each time $t_k$, we will consider the behavior of the process
with respect to the ball $B(o,R_k)$, where the radii $R_k$ are defined by:
\[
R_k = \begin{cases} 0 & k=0, \\ \gamma (m + \sqrt{t_k(K+1-k)}) & k>0,
\end{cases}
\]
where $\gamma>1$ is some constant to be determined during the proof of
\lemref{ub_induction}. Note that $R_k$ is increasing in $k$, and that
$R_K=2\gamma m$.

With the above notations in mind, let $X_k$ (resp.~$Y_k$) be the number of
particles inside (resp.\ outside) $B(o,R_k)$ at time $t_k$. Let $A_m=
A_m(\alpha, \beta,\gamma)$ be the event
\[
A_m = \left\{X_k < \frac{\alpha R_k^d}{t_k} \quad\text{and}\quad Y_k <
  \frac{\beta R_k^d}{t_k} \quad\text{for all $k\leq K$}\right\} \qquad
\text{when $d>2$}.
\]
and
\[
A_m = \left\{X_k < \frac{\alpha R_k^2 \log m}{t_k} \quad\text{and}\quad Y_k
  < \frac{\beta R_k^d \log m}{t_k} \quad\text{for all $k\leq K$}\right\}
\qquad \text{when $d=2$}.
\]

\begin{lemma}\label{L:ub_induction}
  Assume the initial conditions of \thmref{upper_bd}. Then for some choice
  of $\alpha,\beta,\gamma$ we have $\P(A_m) \xrightarrow[m\to\infty]{} 1$.
\end{lemma}

\begin{proof}
  The idea is to inductively bound $X_{k+1},Y_{k+1}$ in terms of $X_k,Y_k$.
  The bound on $X_{k+1}$ is mostly an application of \lemref{low_density}.
  However, to bound $X_{k+1}$ we need to also account for particles that are
  outside the ball $B(o,R_k)$ at time $t_k$, or particles that exit the
  ball $B(o, R_{k+1})$ at some time before time $t_{k+1}$ and re-enter it.
  These quantities can be bounded in terms of $Y_k$ and $X_k$ as well as
  auxiliary quantities introduced soon. A delicate point in the proof is
  that the number of steps of the induction is not fixed ($\ln\ln m$), so
  we make sure that constants do not grow with $m$. Thus all constants
  below depend only on $d$.

  At time $t_0$ our assumptions are that $Y_0=0$ and $X_0\leq C_0 m^d$. For
  the induction step we define two additional quantities: $S_k$ and $Z_k$.
  Let $S_k$ be the number of particles that remain in $B(o,R_k)$ throughout
  the time interval $[t_{k-1},t_k]$. We wish to apply \lemref{low_density}
  to $S_k$. The conditions are clearly satisfied (recall $S_k\leq C_0 m^d <
  C_0 R_k^d$, since $\gamma>1$). Since $t_k-t_{k-1} \geq c t_k$ this will
  imply that with high probability
  \begin{equation}\label{eq:S_bound}
    S_k < C_1 \frac{R_k^d}{t_k}    \qquad\text{ for all } k\le K.
  \end{equation}
  (The probability of failure at each of $\log\log m$ steps is
  exponentially small.)

  Let $Z_k$ be the number of particles located inside $B(o,R_k)$ at time
  $t_k$ that exit $B(o,R_{k+1})$ before time $t_{k+1}$. \lemref{gaussian},
  bounds the escape probability for each of $X_k$ particles inside
  $B(o,R_k)$. Coalescence can only reduce the number of escaping particles,
  so given $X_k$,
  \[
  Z_k \slt \Bin\left( X_k, C_2
    \exp\left(-c_2\frac{(R_{k+1}-R_k)^2}{t_{k+1}-t_k} \right) \right).
  \]
  Here $c_2,C_2$ depend only on $d$, and ``$\slt$'' denotes stochastic
  domination.

  If $k=0$ this implies
  \[
  \E Z_0 \leq C_0 m^d C_2 e^{-c_2 \gamma^2 m^2/t_1} .
  \]
  Since $m^2/t_1 \asymp \ln m$, by making $\gamma$ large enough we obtain
  $\E Z_0 = o(1)$. Then in particular $Z_0=0$ with probability tending to
  $1$. For $k \ge 1$ we use Lemma \ref{bin_dev}, with $\Delta = K^4$ to
  find
  \[
  \P\left(\bigcup_{k=1}^K \left\{ Z_k \ge K^4 + 2\E Z_k \right\} \right)
  \leq K/\sqrt{K^4} = 1/K.
  \]
  Thus with probability at least $1-K^{-1}\to 1$ (as $m\to \infty$), we
  have
  \[
  Z_k < K^4 + 2C_2 X_k \exp\left(-c_2 \frac{(R_{k+1}-R_k)^2}{t_{k+1}-t_k}
  \right) \qquad\text{ for all $k$}.
  \]
  Using $t_{k+1}=e t_k$, $k\geq 1$, together with $\sqrt{t_{k+1}(K-k)} \geq
  \sqrt{e/2} \sqrt{t_k(K+1-k)}$, we conclude that with high probability,
  $Z_0=0$ and
  \begin{equation}\label{eq:ZfromX}
    Z_k < K^4 + 2 C_2 X_k e^{-c_3 \gamma^2(K-k)}
    \qquad \text{for all } k\leq K.
  \end{equation}

  \medskip

  With these preparations in place, we are ready for the induction. Assume
  from here on that $Z_0=0$ and \eqref{eq:S_bound}, \eqref{eq:ZfromX} hold.
  We have, for each $k$, the deterministic bounds
  \begin{align*}
    X_k &\leq S_k + Y_{k-1} + Z_{k-1}, \\
    Y_k &\leq Y_{k-1} + Z_{k-1}.
  \end{align*}
  To see this, note that particles in $B(o,R_k)$ either stayed inside
  ($S$), started outside ($Y$), or exited and returned ($Z$). The bound on
  $Y$ is similar.

  We now carry out an induction over $k$ to bound $X_k,Y_k$ for all
  $k=1,\dots,K$. Suppose that the bound (from the event $A_m$) on $X_j$ and
  $Y_j$ hold for all $j<k$. It follows from \eqref{eq:ZfromX} and the
  inductive hypothesis that
  \begin{align*}
    Y_k \leq \sum_{j<k} Z_j
    &< k K^4 + 2C_2\alpha \sum_{1\leq j<k}
    \frac{R_j^d}{t_j}e^{-c_3\gamma^2(K-j)} \\
    &< k K^4 + 2C_2\alpha \frac{R_k^d}{t_k}
    \sum_{j<k} \frac{t_k}{t_j} e^{-c_3\gamma^2(K-j)} \\
    &= k K^4 + 2C_2\alpha \frac{R_k^d}{t_k} e^{-c_3\gamma^2(K-k)}
    \sum_{j<k} e^{(1-c_3\gamma^2)(k-j)}.
  \end{align*}
  We require $\gamma$ to be large enough that $c_3\gamma^2>2$. Then the
  last sum is at most $1/(e-1)<1$ and so
  \begin{equation}\label{eq:YfromX}
    Y_k < k K^4 + 2C_2 \alpha \frac{R_k^d}{t_k} e^{-c_3\gamma^2(K-k)}.
  \end{equation}
  This proves the induction step for $Y_k$ with any choice of $\beta >
  2C_2\alpha$ (since $k K^4 \leq K^5 \ll R_k^d/t_k$), on the event from
  \eqref{eq:ZfromX}, for all sufficiently large $m$.

  It remains to bound $X_k$, for which we will use the bounds on $S_k$,
  $Y_{k-1}$ and $Z_{k-1}$. We already have
  \[
  Y_{k-1} < (k-1)K^4 + 2C_2 \alpha \frac{R_{k-1}^d}{t_{k-1}}
  e^{-c_3\gamma^2(K+1-k)} < (k-1)K^4 + 2eC_2 \alpha \frac{R_k^d}{t_k}
  e^{-c_3\gamma^2}.
  \]
  Using the induction hypothesis and \eqref{eq:ZfromX} (with $k-1$
  replacing $k$) one finds
  \[
  Z_{k-1} < K^4 + 2C_2 \alpha \frac{R_{k-1}^d}{t_{k-1}}
  e^{-c_3\gamma^2(K+1-k)} < K^4 + 2eC_2 \alpha \frac{R_k^d}{t_k}
  e^{-c_3\gamma^2}.
  \]
  Thus we have
  \begin{align}
    X_k \leq S_k + Y_{k-1} + Z_{k-1} &< C_1 \frac{R_k^d}{t_k} +k K^4
    + 4eC_2 \alpha \frac{R_k^d}{t_k} e^{-c_3\gamma^2} \nonumber \\
    &= k K^4 + \left(C_1 + 4eC_2 \alpha e^{-c_3\gamma^2} \right)
    \frac{R_k^d}{t_k}.
    \label{eq:X_bound}
  \end{align}
  To finish the proof it remains to select $\alpha$ and $\gamma$ (and
  $\beta>2C_2\alpha$) so that
  \[
  C_1 + 4eC_2 \alpha e^{-c_3\gamma^2} < \alpha.
  \]
  This is done by requiring $\gamma$ to satisfy $e^{c_3\gamma^2} > 4eC_2$
  and taking any sufficiently large $\alpha$.

  \medskip

  Turning to the case $d=2$, we proceed along the same lines.
  \lemref{low_density} gives with probability with high probability
  \begin{equation}\label{eq:S_bound2}
    S_k < C_1 \frac{R_k^2 \log m}{t_k}    \qquad\text{ for all } k\le K.
  \end{equation}
  The failure probability is $\log^{-1} m$ at each of $\log\log m$ steps.
  Note that $\log t_i \sim \log m$ for all $i$, so we are not giving much
  away here. Furthermore, w.h.p.\ $Z_0=0$ and \eqref{eq:ZfromX} holds (the
  proof of these facts does not depend on $d$.

  We now repeat the induction. Given \eqref{eq:ZfromX} and the induction
  hypothesis bounds on $X$ we get (as above, with an extra $\log m$ factor)
  \[
  Y_k < k K^4 + 2C_2 \alpha \frac{R_k^2 \log m}{t_k} e^{-c_3\gamma^2(K-k)}.
  \]
  Since the bounds for $S_k$, $Y_{k-1}$ and $Z_{k-1}$ differ from the
  general case only by a $\log m$ on the $R_k^2/t_k$ term, the bound for
  $X_k$ gets the same factor as well.
\end{proof}

\paragraph{Acknowledgements.}
We wish to thank Alan Hammond for useful discussions and careful reading of
a few preliminary drafts. This project was started when all authors were at
the University of British Columbia. Part of it was also carried as the
first author was at the University of Toronto. Further progress was made
during visits to Marseilles through the support of V.L.'s Alfred P. Sloan
fellowship, and to the University of Toronto. The hospitality and support
of these departments is warmly acknowledged.

\bibliographystyle{abbrv}
\bibliography{coalesce}

\begin{minipage}{\textwidth}
  \noindent
  {\bf Omer Angel}: {\tt angel@math.ubc.ca}\\
  Department of Mathematics, University of British Columbia,\\
  Vancouver, BC, V6T 1Z2, Canada

  \bigskip

  \noindent
  {\bf Nathanael Berestycki}: {\tt n.berestycki@statslab.cam.ac.uk} \\
  Statistical Laboratory, DPMMS \\
  University of Cambridge \\
  Wilberforce Rd. Cambridge, CB3 0WB, United Kingdom

  \bigskip

  \noindent
  {\bf Vlada Limic}: {\tt vlada@cmi.univ-mrs.fr}\\
  CNRS - Universit\'e de Provence, Technop\^ole de Ch\^ateau-Gombert \\
  UMR 6632, LATP, CMI \\
  39, rue F.\ Joliot Curie, 13453 Marseille, cedex 13, France
\end{minipage}

\end{document}